# Integration of singular subalgebroids by diffeological groupoids [1]

by Iakovos Androulidakis and Marco Zambon

National and Kapodistrian University of Athens
Department of Mathematics
Panepistimiopolis
GR-15784 Athens, Greece
e-mail: iandroul@math.uoa.gr

KU Leuven
Department of Mathematics
Celestijnenlaan 200B box 2400
BE-3001 Leuven, Belgium.
e-mail: marco.zambon@kuleuven.be

**Abstract**

We establish an integration theory for singular subalgebroids, by diffeological groupoids. To do so, we single out a class of diffeological groupoids satisfying specific properties, and we introduce a differentiation-integration procedure under which they correspond to singular subalgebroids. Our definition of integration distinguishes the holonomy groupoid from the graph, although both differentiate to the original singular subalgebroid: the holonomy groupoid satisfies a certain submersive property, while the graph does not.

## Contents

---

[1] 2020 AMS subject classification: Primary 22A22, 57R55, Secondary 53D17.

# Introduction

Lie algebroids are generalizations of Lie algebras that appear naturally in differential geometry and mathematical physics, and consist essentially of a vector bundle whose space of sections is endowed with a compatible Lie algebra bracket. The sub-objects of Lie algebroids commonly considered are Lie subalgebroids, themselves vector subbundles. *Singular subalgebroids* were introduced in [31] as a more flexible class of sub-objects, and were shown to encompass an abundance of examples.

## Subject of the paper

Our primary purpose in this paper is to develop a theory of *integration* for singular subalgebroids. Here "integration" is taken in the same sense as a Lie group is an integration of a Lie algebra. To this end, first let us recall the definition of singular subalgebroid, as well as the constructions given in [31].

**Definition.** Let $A$ be a Lie algebroid over the manifold $M$. A **singular subalgebroid** is a $C^\infty(M)$-submodule $\mathcal{B}$ of compactly supported sections of $A$, which is locally finitely generated and involutive with respect to the Lie algebroid bracket.

Notice that the notion of singular subalgebroid is a relative one, since it depends on the ambient Lie algebroid $A$.

In fact, the notion of *singular foliation* $(M, \mathcal{F})$ from [5] is our guiding example of singular subalgebroid, where the ambient Lie algebroid is the tangent bundle $TM$. Recall that $TM$ can be integrated by more than one Lie groupoid, for instance the fundamental groupoid and the pair groupoid. In [31] it was noticed that the holonomy groupoid $H(\mathcal{F})$ constructed in [5] really takes into account one more choice; that is, the choice of Lie groupoid integrating $TM$. Indeed, to every triple $(\mathcal{B}, A, G)$ consisting of a singular subalgebroid $\mathcal{B}$ of a Lie algebroid $A$ and a choice of Lie groupoid $G$ integrating the latter, [31] attached a topological groupoid $H^G(\mathcal{B})$. In the case of a singular foliation, applying this construction to $(\mathcal{F}, TM, M \times M)$ recovers $H(\mathcal{F})$. In fact, we call $H^G(\mathcal{B})$ itself *holonomy groupoid*, because it is a kind of pullback of a holonomy groupoid for a singular foliation on $G$.

To get an idea of what the integration of singular subalgebroids entails, we look at two cases. They provide an indication that the integration problem for singular subalgebroids is even more relative than their definition: one not only needs to prescribe the ambient Lie algebroid $A$, but also a Lie groupoid $G$ integrating $A$.

- Consider a wide Lie subalgebroid $B$ of an integrable Lie algebroid $A$. One can consider $B$ itself as a Lie algebroid, with the structure induced by $A$. However doing so one loses information, while we wish to remember that $B$ is a sub-object of the ambient Lie algebroid $A$. In this context, Moerdijk-Mrčun fixed a Lie groupoid $G$ integrating $A$, and in [24, Thm. 2.3] they classified the Lie groupoids[2] $H$ which admit a Lie groupoid morphism[3] $H \to G$ integrating the inclusion $B \hookrightarrow A$. They called them *integrals of $B$ over $G$*. Further they showed that there is a minimal such Lie groupoid $H$, which agrees with the holonomy groupoid $H^G(B)$. We emphasize again that the focus is not on $B$ itself as a Lie algebroid, but as a sub-object of $A$.
- Given a triple triple $(\mathcal{B}, A, G)$ as above, the holonomy groupoid $H^G(\mathcal{B})$ comes together with a canonical morphism of topological groupoids $\Phi : H^G(\mathcal{B}) \to G$ [31].

It turns out that an essential aspect of an integration theory for singular subalgebroids is the specification of a correct notion of *smoothness* for the integrating groupoids. In the approach we propose, which we outline below, the integrals of singular subalgebroids are diffeological groupoids [16, §8.3], i.e. groupoids endowed with a compatible diffeology. This was already announced in [31], and diffeology in the context of the holonomy groupoid was later used by other authors [14] [21] [22]. In full generality, diffeological groupoids are objects that are too flexible in order to carry out differential geometric constructions.

Therefore, one of the main challenges of this paper is to identify conditions on diffeological groupoids which, once imposed, are sufficient to carry out the integration theory of singular subalgebroids, while not being too restrictive. We present these conditions in §2.3, and we give an overview at the end of this introduction.

**Remark.** Diffeological groupoids occur in several contexts in the literature, including general relativity [11], quantization [3][4] and (higher) Lie theory [12][26][10][2]. We emphasize that in this paper we only consider diffeological groupoids whose space of objects is a smooth manifold, and which come together with a morphism to a Lie groupoid. This simplifies many of our considerations and allows us to avoid some technical aspects of the theory of diffeological spaces (for instance, we never have to address their tangent spaces).

[2] $H$ necessarily integrates $B$ as a Lie algebroid.
[3] The morphisms $H \to G$ are not injective in general.

## Results

We explain our approach to the integration question. Let us fix a triple $(\mathcal{B}, A, G)$ consisting of a singular subalgebroid $\mathcal{B}$ of a Lie algebroid $A$ and a Lie groupoid $G$ integrating the latter.

The first step is to get a good grasp on *differentiation.* That is because a meticulous examination of differentiation will reveal the diffeological properties which are necessary for our integration theory. The holonomy groupoid $H^G(\mathcal{B})$ provides a guide. Hence we ask:

*I) Can we find a differentiation procedure which recovers $\mathcal{B}$ from $H^G(\mathcal{B})$?*

The "link" between the holonomy groupoid and the Lie algebroid $A$ is provided by the canonical morphism $\Phi\colon H^G(\mathcal{B}) \to G$. We will show that $\mathcal{B}$ is recovered as the (compactly supported) elements of

$$\left\{ \frac{d}{d\lambda}|_{\lambda=0}(\Phi \circ \mathbf{b}_\lambda) \colon \ \{\mathbf{b}_\lambda\}_{\lambda\in I} \text{ family of smooth global bisections for } H^G(\mathcal{B}) \text{ s.t. } \mathbf{b}_0 = Id_M \right\}.$$

An important point is to make sense of what it means for a bisection to be "smooth". This is done using the natural diffeological structure on $H^G(\mathcal{B})$, arising from the fact that $H^G(\mathcal{B})$ is a quotient of a disjoint union of manifolds.

**Remark.** We point out that we are using two well-known and fundamental ideas about differentiation.

- First, a global integration question was considered by Souriau: Given a compact manifold $M$, in what sense can we apply the Lie functor to the group of diffeomorphisms $Diff(M)$ in order to obtain the Lie algebra $\mathfrak{X}(M)$ of vector fields? Souriau achieved this differentiation by introducing the notion of *diffeology* in [27].
- The second key to the differentiation of $H^G(\mathcal{B})$ comes from H. Abels [1, Appendix], who showed that in order to recover the infinitesimal information (the Lie algebra) one just needs to make good sense of 1-parameter families, and this can be done even when we work with just a topological group. In fact, this idea can be traced back to Palais.

In §3, guided by the above example of the holonomy groupoid, we single out a class of objects that give rise to singular subalgebroids by a differentiation procedure as the one above. They are given by pairs $(H, \Psi)$ where $H$ is a diffeological groupoid whose base space is a smooth manifold and $\Psi\colon H \to G$ a morphism covering $Id_M$, satisfying certain properties[4] modelled on those of the holonomy groupoid. Our approach is minimalistic, in the sense that the properties we require are the minimal assumptions under which we were able to extract a well-defined singular subalgebroid from our data. This leads to the definition of differentiation (Definition 3.3). As expected, $H^G(\mathcal{B})$ differentiates to $\mathcal{B}$ in this sense.

With the notion of differentiation at hand, we can finally address the integration question:

*II) What is the correct notion of integration for a singular subalgebroid $\mathcal{B}$?*

Our answer is more refined than just pairs $(H, \Psi)$ with the expected requirement that they differentiate to $\mathcal{B}$: we impose two more conditions that make the notion of integration as well-behaved as possible, while retaining all the desired classes of examples. We paraphrase Definition 4.1.

[4]The "holonomy-like" property of Def. 2.12), and Assumption 3.1.

**Definition.** Fix a triple $(\mathcal{B}, A, G)$ consisting of a singular subalgebroid $\mathcal{B}$ of a Lie algebroid $A$ and a Lie groupoid $G$ for $A$. Let $H$ be a diffeological groupoid over $M$, and $\Psi\colon H \to G$ a smooth morphism of diffeological groupoids covering $Id_M$. We say that $(H, \Psi)$ is an **integral of $\mathcal{B}$ over $G$** if:

a) $(H, \Psi)$ differentiates to $\mathcal{B}$ (Def. 3.3).

b) The diffeology of $H$ is generated by open maps (Def. 2.18).

c) The morphism $\Psi$ is *coarsely injective.*

Condition b) can be viewed as a property on the topology of $H$, and condition c) is there to exclude groupoids that are "larger than necessary". Notice that we adopt the terminology "integral", to distinguish this notion from the usual integration by Lie groupoids. To test the validity of our notion of integral, we apply it to the main classes of examples. We summarize the results as follows (Thm. 4.19, Prop. 4.7, Prop. 4.17, Rem. 4.18). Notice that items ii), iii), iv) are in increasing order of generality. Here we use the term "smooth integral" to denote integrals that are Lie groupoids.

**Theorem.** The following are classes of examples of integrals.

i) For any singular subalgebroid $\mathcal{B}$, the holonomy groupoid $H^G(\mathcal{B})$ is an integral.

ii) When $A$ is a Lie algebra (thus $\mathcal{B} =: B$ is a Lie subalgebra), the integrals of $B$ are all smooth. Indeed, the integrals of $B$ over $G$ are exactly the Lie groups covering the connected Lie subgroup of $G$ integrating the Lie subalgebra $B$.

iii) When $\mathcal{B} = \Gamma_c(B)$ for a wide Lie subalgebroid $B$, the smooth integrals coincide with the integrals of $B$ over $G$ in the sense of Moerdijk-Mrčun [24].

iv) A singular subalgebroid $\mathcal{B}$ admits a smooth integral iff it is a projective, i.e. $\mathcal{B} \cong \Gamma_c(B)$ for a Lie algebroid $B$. The pair $(H^G(\mathcal{B}), \Phi)$ is a smooth integral of $\mathcal{B}$ over $G$, and further it is the minimal smooth integral: the smooth integrals are exactly the coverings of the holonomy groupoid.

For any singular subalgebroid $\mathcal{B}$, one can also consider the *graph* of $\mathcal{B}$, a set-theoretic subgroupoid of $G$. It is defined analogously to the graph of a foliation, which is the collection of pairs of points belonging to the same leaf. The graph is quite a coarse object, but as a quotient of the holonomy groupoid $H^G(\mathcal{B})$ it carries a diffeological groupoid structure, which turns out to differentiate to $\mathcal{B}$. Since the inclusion of the graph in $G$ is an injective map, the graph is the smallest diffeological groupoid differentiating to $\mathcal{B}$. We show in §5 that the openness property b) is *not* satisfied by the graph in general, hence the graph is not an integral in general. This supports the validity of our theory of integration: it allows to distinguish the holonomy groupoid, which is always an integral, from the graph of the singular subalgebroid, which typically fails to be an integral.

**Remark.** For the holonomy groupoid, the openness property in the above Definition is implied by a kind of submersive property for the plots of the diffeological structure, which we call "local subduction property". We define it in §2.3.3, and view it as a smoothness property of the diffeological structure. The investigation of a submersive property was suggested to us by Georges Skandalis. It has independent interest, in view of the construction of operator algebras from diffeological groupoids, which are interesting in Noncommutative Geometry. However we do not address this here.

**Remark.** In this paper we introduce several properties of two kinds of mathematical objects:

– a diffeological groupoid $H \rightrightarrows M$,

– a smooth morphism of diffeological groupoids $\Psi : H \to G$ covering $Id_M$, where $G$ is a Lie groupoid over $M$.

This table gives an overview of the main properties. The properties in the first line are needed to make sense of *differentiation*, while all four are required in our definition of *integral.*

| *Properties of H* | *Properties of* $\Psi\colon H \to G$ |
|---|---|
| holonomy-like (Def. 2.12) | Assumption 3.1 on 1-parameter groups |
| open map property (Def. 2.18) | $\Psi$ is coarsely injective |

Further, the open map property is related as follows to other properties of diffeological groupoids we introduce in §2.3:

local subduction property $\Rightarrow$ open map property $\Rightarrow$ source-submersive property.

In Appendix B we show that whenever $\Psi$ is coarsely injective one has:

Assumption 3.1 on 1-parameter groups $\Leftrightarrow$ Assumption B.3 on liftability.

**Structure of the paper:** The article starts in §1 by recalling singular subalgebroids and their associated holonomy groupoids from [31]. In §2 we define diffeological groupoids and discuss the properties which are relevant to the integration of singular subalgebroids. §3 gives the differentiation process for diffeological groupoids endowed with a morphism to a Lie groupoid (most proofs are collected in Appendix A). The integration is discussed in §4. Finally, the graph of a singular subalgebroid is discussed in §5.

**Relation to the literature:** Given a singular subalgebroid $\mathcal{B}$ of a Lie algebroid $A$ and a Lie groupoid $G$ integrating the latter, as explained above, we carry out a theory of integration of $\mathcal{B}$ that keeps track of the fact that $\mathcal{B}$ is a sub-object of $A$, in analogy with [24]. One could instead focus on the intrinsic structure carried by $\mathcal{B}$, namely the structure of Lie-Rinehart algebra inherited from $A$. Its integration is addressed in a special case in [14], and in generality by Villatoro in [28]; in this context, however, it is not known whether the resulting diffeological groupoids allow to recover $\mathcal{B}$ by some sort of differentiation procedure.

**Conventions and notation:** All Lie groupoids are assumed to be source connected. Given a Lie groupoid $G \rightrightarrows M$, we denote by $\mathbf{t}$ and $\mathbf{s}$ its target and source maps, and by $\iota\colon G \to G$ the inversion map. We denote by $1_x \in G$ the identity element corresponding to a point $x \in M$, and by $1_M \subset G$ the submanifold of identity elements. Two elements $g, h \in G$ are composable if $\mathbf{s}(g) = \mathbf{t}(h)$. We use the term bisection to denote a right-inverse $\mathbf{b}$ to $\mathbf{s}$ defined on an open subset of $M$, such that $\mathbf{t} \circ \mathbf{b}$ is a diffeomorphism onto its image. We identify the Lie algebroid of $G$, which we denote by $AG$, with $\ker(d\mathbf{s})|_M$.

**Acknowledgements:** I.A. thanks Georges Skandalis for various discussions and suggestions. He was paritally supported by COST CaLISTa CA21109. He would also like to thank Collegium Helveticum and the Institute of Mathematics of the University of Zurich, for their hospitality during the revision of this article. M.Z. thanks Dan Christensen, Claire Debord, Alfonso Garmendia and Joel Villatoro for explanations related to this work. We thank the anonymous referee for pointing out the preprint [9]. I.A. was partially supported by a Marie Curie Career Integration Grant PCI09-GA-2011-290823 (Athens), SCHI 525/12-1 (DFG/Germany), thanks the Collegium Helveticum, and acknowledges partial support of the SNF Grant No. 200020 192080 (Switzerland) and by the NCCR SwissMAP, funded by the Swiss SNF. M.Z. was partially supported by grants MTM2011-22612 and ICMAT Severo Ochoa SEV-2011-0087 (Spain), Pesquisador Visitante Especial grant 88881.030367/2013-01 (CAPES/Brazil), IAP Dygest, Methusalem grants METH/15/026 and METH/21/03 - long term structural funding of the Flemish Government, the FNRS and FWO under EOS projects G0H4518N and G0I2222N, the FWO research projects G083118N and G0B3523N (Belgium).

# 1 Background on singular subalgebroids, holonomy groupoids and diffeology

In this section we briefly recall from [31] (see also [8]) the notion of a singular subalgebroid in §1.1, and the associated holonomy groupoid in §1.2-1.3. Based on this, in §2, we will give the description of the holonomy groupoid as a diffeological space, which is crucial for our later results. In order to make the exposition of these results self-contained, in §1.4 we recall diffeological spaces from [16] and [27].

Let $A$ be a Lie algebroid over a smooth manifold $M$. Throughout the sequel we assume $A$ is integrable and choose a (source-connected) Lie groupoid $G \rightrightarrows M$ integrating $A$. We will assume that $G$ is Hausdorff.

## 1.1 Singular subalgebroids and examples

**Definition 1.1.** [31, Def. 1.1] A **singular subalgebroid** of $A$ is an involutive, locally finitely generated $C^\infty(M)$-submodule $\mathcal{B}$ of $\Gamma_c(A)$, the compactly supported sections of $A$.

We briefly discuss the main examples of singular subalgebroids.

**Examples 1.2** (Singular subalgebroids).

a) Let $M$ be a smooth manifold and put $A = TM$. A singular subalgebroid of $TM$ is a **singular foliation** $(M, \mathcal{F})$ in the sense of [5].

b) Let $\psi\colon E \to A$ be a morphism of Lie algebroids covering the identity. It induces a morphism of $C^\infty(M)$-modules between the respective modules of compactly supported sections, which we still denote $\psi$. The image $\mathcal{B} := \psi(\Gamma_c(E))$ is a singular subalgebroid of $A$. We say that $\mathcal{B}$ **arises from the Lie algebroid morphism** $\psi$.

c) Special cases of the previous item are the **projective** singular subalgebroids $\mathcal{B}$ of $A$, namely those for which there exists a vector bundle $B$ over $M$ such that $\Gamma_c(B) \cong \mathcal{B}$ as $C^\infty(M)$-modules. There is a Lie algebroid structure on $B$ and an almost injective Lie algebroid morphism $\tau\colon B \to A$ inducing the isomorphism $\Gamma_c(B) \cong \mathcal{B}$, and this data is unique, see [6]. In particular, $\mathcal{B}$ arises from the Lie algebroid morphism $\tau$.

d) Specializing further, we can consider $\mathcal{B} = \Gamma_c(B)$ where $B$ is a **wide Lie subalgebroid** of $A$ (that is, a vector subbundle $B \to M$ whose sections are closed with respect to the Lie bracket [23, Def. 3.3.21]).

e) Let $N$ be a closed embedded submanifold of $M$ and $B \to N$ a **Lie subalgebroid** of $A$ ([23, Def. 4.3.14]). Then

$$\mathcal{B} := \{\boldsymbol{\alpha} \in \Gamma_c(A) : \boldsymbol{\alpha}|_N \subset B\}$$

is a singular subalgebroid of $A$. Let $n = dim(N)$ and $b = rank(B)$. It was shown in [31, §1.3.3] that when $N$ is a hypersurface, $\mathcal{B}$ is projective. On the other hand, if $codim(N) \geqslant 2$ and $B \neq A|_N$, then $\mathcal{B}$ is not projective.

Let us recall from [31, §1] the following singular foliations associated with a singular subalgebroid $\mathcal{B}$:

a) *The singular foliation* $\overrightarrow{\mathcal{B}}$ *on* $G$*:* Since the Lie algebroid $A$ integrates to the Lie groupoid $G$, every section $\boldsymbol{\alpha}$ of $A$ corresponds to a right-invariant vector field $\overrightarrow{\boldsymbol{\alpha}}$ of $G$ (see [23, §3.5]). Whence the $C^\infty(M)$-module $\mathcal{B}$ gives rise to the involutive $C^\infty_c(G)$-module $\overrightarrow{\mathcal{B}}$ generated by $\{\overrightarrow{\boldsymbol{\alpha}} : \boldsymbol{\alpha} \in \mathcal{B}\}$. Notice that $\overrightarrow{\mathcal{B}}$ consists of compactly supported vector fields [5] That is, $(G, \overrightarrow{\mathcal{B}})$ is a singular foliation in the sense of [5].

b) Likewise, every $\boldsymbol{\alpha} \in \mathcal{B}$ defines a left-invariant vector field $\overleftarrow{\boldsymbol{\alpha}}$ of $G$ and a singular foliation $\overleftarrow{\mathcal{B}}$.

c) *The singular foliation* $\mathcal{F}_{\mathcal{B}}$ *on* $M$*:* Let $\rho : A \to TM$ be the anchor map of $A$. Since $\mathcal{B}$ is locally finitely generated and involutive, so is the $C^\infty(M)$-module $\mathcal{F}_{\mathcal{B}}$ generated by $\{\rho(\boldsymbol{\alpha}) : \boldsymbol{\alpha} \in \mathcal{B}\}$; whence $\mathcal{F}_{\mathcal{B}}$ is a singular foliation of $M$. Note that $\rho(\boldsymbol{\alpha})$ is compactly supported. Recall that $\overrightarrow{\boldsymbol{\alpha}}$ is $\mathbf{t}$-related with $\rho(\boldsymbol{\alpha})$, therefore the vector field $\overrightarrow{\boldsymbol{\alpha}}$ is complete.

Last, we will need the *fibers* of $\mathcal{B}$: For every $x \in M$, put $I_x = \{f \in C^\infty(M) : f(x) = 0\}$. Then $I_x\mathcal{B}$ is an ideal of $\mathcal{B}$. Since $\mathcal{B}$ is locally finitely generated, the quotient $\mathcal{B}_x = \mathcal{B}/I_x\mathcal{B}$ is a finite-dimensional vector space. We call $\mathcal{B}_x$ the **fiber of $\mathcal{B}$ at** $x \in M$.

## 1.2 Bisubmersions

The main tool for the construction of the holonomy groupoid of a singular subalgebroid is the notion of bisubmersion. We recall from [31] this notion here, as well as some results that we'll need in this sequel.

Let $U, V$ be manifolds and $\varphi : U \to V$ a smooth map. Let $\mathcal{E}$ be a $C^\infty(V)$-submodule of $\mathfrak{X}_c(V)$. Denote $\varphi^! TV$ the vector bundle on $U$ obtained as the the pullback of the tangent bundle. Put:

a) $\varphi^*(\mathcal{E})$ the $C^\infty(U)$-submodule of $\Gamma_c(U, \varphi^! TV)$ generated by $f(\xi \circ \varphi)$ with $f \in C^\infty_c(U)$ and $\xi \in \mathcal{E}$.

b) $\varphi^{-1}(\mathcal{E})$ the $C^\infty(U)$-submodule $\{X \in \mathfrak{X}_c(U) : d\varphi(X) \in \varphi^*(\mathcal{E})\}$ of $\mathfrak{X}_c(U)$.

**Definition 1.3.** [31, Prop. 2.4] Let $\mathcal{B}$ be a singular subalgebroid of the Lie algebroid $A$, as in Definition 1.1. A **bisubmersion** for $\mathcal{B}$ is a triple $(U, \varphi, G)$ where $U$ is a manifold and $\varphi\colon U \to G$ a smooth map such that

[5]This is the reason why we assume the Lie groupoid $G$ to be Hausdorff: the compact support condition on non-Hausdorff manifolds leads to somewhat strange behaviours. For instance, compact sets might not be closed. We assume $G$ is Hausdorff, in order to avoid complicating the exposition. (Allowing $G$ to be non-Hausdorff, one would have to work with sheaves of compactly supported sections on open subsets of $G$).

i) $\mathbf{s}_U := \mathbf{s} \circ \varphi$ and $\mathbf{t}_U := \mathbf{t} \circ \varphi \colon U \to M$ are submersions,

ii) for every $\boldsymbol{\alpha} \in \mathcal{B}$, there is $Z \in \mathfrak{X}(U)$ which is $\varphi$-related to $\overrightarrow{\boldsymbol{\alpha}}$ and $W \in \mathfrak{X}(U)$ $\varphi$-related to $\overleftarrow{\boldsymbol{\alpha}}$,

iii) $\varphi^{-1}(\overrightarrow{\mathcal{B}}) = \Gamma_c(U, \ker d\mathbf{s}_U)$ and $\varphi^{-1}(\overleftarrow{\mathcal{B}}) = \Gamma_c(U, \ker d\mathbf{t}_U)$.

**Definition 1.4.** Let $(U_i, \varphi_i, G)$, $i = 1, 2$ two bisubmersions for $\mathcal{B}$. A **morphism** of bisubmersions is a smooth map $f : U_1 \to U_2$ such that $\varphi_2 \circ f = \varphi_1$.

We now recall several constructions from [31, §2.4, 2.5], which give rise to bisubmersions and morphisms thereof.

**Examples 1.5** (Bisubmersions and their morphisms)**.** a) Recall that, given any $\boldsymbol{\alpha} \in \mathcal{B}$ the vector field $\overrightarrow{\boldsymbol{\alpha}}$ on $G$ is complete. Let $x \in M$ and $\boldsymbol{\alpha}_1, \ldots, \boldsymbol{\alpha}_n \in \mathcal{B}$ such that $[\boldsymbol{\alpha}_1], \ldots, [\boldsymbol{\alpha}_n]$ span the fiber $\mathcal{B}_x$. The associated **path-holonomy bisubmersion** is the map

$$\varphi \colon U_0 \to G, (\lambda, x) \mapsto \exp_x \sum \lambda_i \overrightarrow{\boldsymbol{\alpha}_i},$$

where $U_0$ is a neighborhood of $(0, x)$ in $\mathbb{R}^n \times M$ in which the map $\mathbf{t} \circ \varphi$ is a submersion, and exp denotes the time-1 flow.

b) In the previous item, we say that $(U_0, \varphi, G)$ is **minimal** if $[\boldsymbol{\alpha}_1], \ldots, [\boldsymbol{\alpha}_n]$ are a basis of $\mathcal{B}_x$.

c) The **inverse** of $(U, \varphi, G)$ is the bisubmersion $\iota \circ \varphi : U \to G$, where $\iota : G \to G$ is the inversion map of $G$.

d) Let $(U_i, \varphi_i, G)$, $i = 1, 2$ be two bisubmersions for $\mathcal{B}$. Put $\mathbf{s}_{U_i} = \mathbf{s} \circ \varphi_i$ and $\mathbf{t}_{U_i} = \mathbf{t} \circ \varphi_i$. The **composition** of $(U_1, \varphi_1, G)$ and $(U_2, \varphi_2, G)$ is the bisubmersion $m \circ (\varphi_1, \varphi_2) : U_1 \times_{\mathbf{s}_{U_1}, \mathbf{t}_{U_2}} U_2 \to G$, where $m$ is the multiplication map of $G$. This bisubmersion is denoted $U_1 \circ U_2$.

e) Let $(U, \varphi, G)$ a path-holonomy bisubmersion for $\mathcal{B}$. The map $\kappa : U \to U, \kappa(\lambda, x) = (-\lambda, \mathbf{t}_U(\lambda, x))$ is a diffeomorphism and the following diagram commutes:

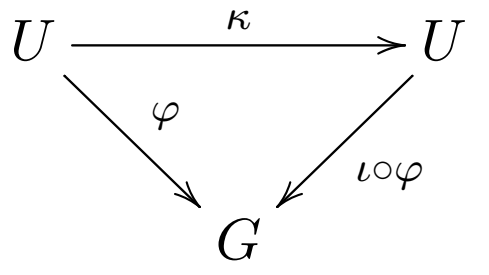

(1.1)

(In particular, $\mathbf{s}_U \circ \kappa = \mathbf{t}_U$ and $\mathbf{t}_U \circ \kappa = \mathbf{s}_U$.) In other words, $\kappa$ is an isomorphism of bisubmersions between $(U, \varphi, G)$ and its inverse bisubmersion. The relevance of this fact is that it simplifies the definition of *path-holonomy atlas* we present below.

There is also the notion of atlas of bisubmersions. Here we define it only for path-holonomy bisubmersions; the general definition is given in [31, Appendix C].

**Definition 1.6.** a) Consider a family $(U_i, \varphi_i, G)_{i \in I}$ of minimal path-holonomy bisubmersions such that $M = \cup_{i \in I} \mathbf{s}_i(U_i)$. Let $\mathcal{U}^G$ be the collection of all such bisubmersions, together with their inverses and finite compositions. We call $\mathcal{U}^G$ a **path-holonomy atlas** of $\mathcal{B}$ associated with $G$.

b) Let $(U, \varphi, G)$ a bisubmersion and $u \in U$. We say that $(U, \varphi, G)$ is **adapted** to $\mathcal{U}^G$ at $u$ if there exists $i \in I$, an open subset $U' \subseteq U$ containing $u$ and a morphism of bisubmersions $U' \to U_i$.

c) We say that $(U, \varphi, G)$ is adapted to $\mathcal{U}^G$ if every element $u \in U$ is adapted to $\mathcal{U}^G$.

### 1.2.1 Bisections

The notion of bisection [31, Def. 2.30] of a bisubmersion will be crucial in the sequel.

**Definition 1.7.** a) A **bisection** of $(U, \varphi, G)$ is a locally closed submanifold $\mathbf{u}$ of $U$ such that the restrictions of both $\mathbf{s}_U$ and $\mathbf{t}_U$ to $\mathbf{u}$ are diffeomorphisms from $\mathbf{u}$ onto open subsets of $M$. Then $\varphi(\mathbf{u})$ is a bisection of $G$ and the map $\varphi : \mathbf{u} \to \varphi(\mathbf{u})$ is a diffeomorphism

b) Let $u \in U$ and $\mathbf{c}$ a bisection of $G$. We say that $\mathbf{c}$ is **carried** by $(U, \varphi, G)$ at $u$ if there exists a bisection $\mathbf{u}$ of $U$ such that $u \in \mathbf{u}$ and $\varphi(\mathbf{u})$ is an open subset of $\mathbf{c}$.

Let $(U, \varphi, G)$ and $(U_i, \varphi_i, G)$, $i = 1, 2$ be bisubmersions. The next properties were proven in [31, §2.5] and [31, §3.1]:

a) Let $u \in U$ and $\mathbf{c}$ a bisection of $G$ carried by $(U, \varphi, G)$ at $u$. Then $\mathbf{c}^{-1}$ is carried by the inverse bisubmersion $(U, \iota \circ \varphi, G)$ at $u$.

b) Let $u_i \in U_i$, $i = 1, 2$ be such that $\mathbf{s}_{U_1}(u_1) = \mathbf{t}_{U_2}(u_2)$ and let $\mathbf{c}_i$ be bisections of $G$ carried by $(U_i, \varphi_i, G)$ at $u_i$ respectively. Then $\mathbf{c}_1 \cdot \mathbf{c}_2$ is carried by the composition $(U_1 \circ U_2, \varphi_1 \cdot \varphi_2, G)$ at $(u_1, u_2)$.

c) Let $(U_0, \varphi_0, G)$ be a minimal path-holonomy bisubmersion. Let $(U, \varphi, G)$ be any bisubmersion carrying the identity bisection at $u \in U$. Then there exists an open neighborhood $U'$ of $u$ in $U$ and a submersion $g : U' \to U_0$ which is a morphism of bisubmersions.

The next properties [31, Cor. 3.3] are crucial for the construction of the holonomy groupoid.
Let $(U_i, \varphi_i, G)$, $i = 1, 2$ be bisubmersions of $\mathcal{B}$ and $u_i \in U_i$ such that $\varphi_1(u_1) = \varphi_2(u_2)$.

a) If the identity bisection $1_M$ is carried by $U_i$ at $u_i$, for $i = 1, 2$, there exists an open neighborhood $U_1'$ of $u_1$ in $U_1$ and a morphism of bisubmersions $f \colon U_1' \to U_2$ such that $f(u_1) = u_2$.

b) If there is a bisection of $G$ carried by both $U_1$ at $u_1$ and by $U_2$ at $u_2$, there exists an open neighborhood $U_1'$ of $u_1$ in $U_1$ and a morphism of bisubmersions $f \colon U_1' \to U_2$ such that $f(u_1) = u_2$.

c) If there is a morphism of bisubmersions $g \colon U_1 \to U_2$ such that $g(u_1) = u_2$, then there exists an open neighborhood $U_2'$ of $u_2$ and a morphism of bisubmersions $f \colon U_2' \to U_1$ such that $f(u_2) = u_1$.

## 1.3 Construction of the holonomy groupoid and examples

Analogously to the case of singular foliations [5], the holonomy groupoid of a singular subalgebroid $\mathcal{B}$ is a quotient. Its construction starts from the following observations:

- The last three properties above show hat for $u_1 \in (U_1, \varphi_1, G), u_2 \in (U_2, \varphi_2, G)$ the relation

$$u_1 \sim u_2 \Leftrightarrow \text{there is an open neighborhood } U_1' \text{ of } u_1, \text{ there is a morphism of bisubmersions } f \colon U_1' \to U_2 \text{ such that } f(u_1) = u_2$$

is an equivalence relation.

- This equivalence relation can be also expressed as follows [31, Remark 3.6]:

$$u_1 \sim u_2 \Leftrightarrow \varphi_1(u_1) = \varphi_2(u_2), \exists \text{ bisections } \mathbf{u}_i \text{ through } u_i \text{ s.t. } \varphi_1(\mathbf{u}_1) = \varphi_2(\mathbf{u}_2).$$

This allows us to give the following definition.

**Definition 1.8.** [31, Def. 3.5] Let $\mathcal{B}$ a singular subalgebroid of $A = AG$ and $\mathcal{U}^G$ an associated path-holonomy atlas. The **holonomy groupoid of $\mathcal{B}$ over $G$** is

$$H^G(\mathcal{B}) := \coprod_{U \in \mathcal{U}^G} U/\sim$$

The holonomy groupoid $H^G(\mathcal{B})$ has the following algebraic structure. Denote $\natural\colon \coprod_{U \in \mathcal{U}^G} U \to H^G(\mathcal{B})$ the quotient map, and $q_i := \natural|_{U_i}$. Endowed with the quotient topology, $H^G(\mathcal{B})$ carries a topological groupoid structure with space of objects $M$. Let us describe its algebraic structure.

a) The source and target maps $\mathbf{s}_H, \mathbf{t}_H : H^G(\mathcal{B}) \to M$ are determined by $\mathbf{s}_H \circ q_i = \mathbf{s}_i$ and $\mathbf{t}_H \circ q_i = \mathbf{t}_i$ for all $i \in I$. The multiplication is determined using the composition of bisubmersions, as follows: $q_i(u)q_j(v) = q_{U_i \circ U_j}(u, v)$ for all $i, j$.

b) There is a canonical morphism of topological groupoids $\Phi\colon H^G(\mathcal{B}) \to G$, $\Phi([u]) = \varphi(u)$, where $u$ is any point in a bisubmersion $(U, \varphi, G)$ belonging to the path-holonomy atlas $\mathcal{U}^G$.

c) For every bisubmersion $(U, \varphi, G)$ adapted to $\mathcal{U}^G$ there is a map $q_U : U \to H^G(\mathcal{B})$ such that, for every local morphism of bisubmersions $f : U' \subseteq U \to U_i$ and every $u \in U'$ we have $q_U(u) = q_i(f(u))$.

One can show that the topological groupoid $H^G(\mathcal{B})$ does not depend on the choice of path-holonomy atlas [31, Remark C.4].

Usually the topology of the holonomy groupoid is quite bad, however it is always longitudinally smooth (*cf.* [6]). In §3 we will show that it is also a diffeological groupoid. For the moment, let us show the following topological property.

**Proposition 1.9.** The quotient map $\natural : \coprod_{U \in \mathcal{U}^G} U \to H^G(\mathcal{B})$ is an open map.

Later we will recover this proposition in an alternative way, see Remark 1.17 b) and Proposition 2.24. Although we will not need it in the sequel, let us point out that proposition 1.9 is valid for any atlas. Indeed, the fact that $\mathcal{U}$ is a path holonomy atlas is not used in the proof below.

*Proof.* Fix an open subset $A \subseteq \coprod_{U \in \mathcal{U}^G} U$. We need to prove that the saturation $\natural^{-1}(\natural(A)) \subseteq \coprod_{U \in \mathcal{U}^G} U$ is open as well. Recall that $\natural^{-1}(\natural(A)) = \{u \in \coprod_{U \in \mathcal{U}^G} U : \natural(u) = \natural(a) \text{ for some } a \in A\}$. Pick an element $u \in \natural^{-1}(\natural(A))$. We need to show that there is an open neighborhood $B$ of $u$ in $\coprod_{U \in \mathcal{U}^G} U$ such that $B \subset \natural^{-1}(\natural(A))$. Pick $a \in A$ such that $\natural(u) = \natural(a)$ and assume that $u, a$ belong to path-holonomy bisubmersions $U_u, U_a$ respectively. Since $\natural(u) = \natural(a)$, by definition of the equivalence relation in the path-holonomy atlas, there is an open subset $\widetilde{U}_u \subseteq U_u$ containing the element $u$, and a morphism of bisubmersions $f : \widetilde{U}_u \to U_a$ such that $f(u) = a$. Since $A$ is open, we can arrange that $f(\widetilde{U}_u) \subset A$ by shrinking the domain if necessary. Then, every other element $u' \in \widetilde{U}_u$ satisfies $\natural(u') = \natural(f(u'))$, and since $f(u') \in A$ it follows that $u' \in \natural^{-1}(\natural(A))$. So the open subset $\widetilde{U}_u$ of $\coprod_{U \in \mathcal{U}^G} U$ is contained in $\natural^{-1}(\natural(A))$. □

**Examples 1.10** (Holonomy groupoids)**.** The following examples of holonomy groupoids were given in [31], but we'll use them in this sequel as well. Fix a source-connected Lie groupoid $G$ integrating a Lie algebroid $A$ (so $A = AG$).

a) Let $\phi : K \to G$ be a morphism of (source-connected) Lie groupoids covering the identity on $M$. Then $\phi : K \to G$ is a bisubmersion for the singular subalgebroid $\mathcal{B} := \phi_*(\Gamma_c(AK))$ of $AG$. (This is [31, Prop 2.13].) In [31, Prop 3.12] it is shown that in this case we have:

   i) $H^G(\mathcal{B}) = K/\mathcal{I}$ as topological groupoids, where

   $$\mathcal{I} := \{k \in K : \exists \text{ a (local) bisection } \mathbf{u} \text{ through } k \text{ such that } \phi(\mathbf{u}) \subset 1_M\}$$

   ii) the canonical map $\Phi\colon H^G(\mathcal{B}) \to G$ coincides with the map $K/\mathcal{I} \to G$ induced by $\phi$.

b) A special case of the above item is the following: Let $B$ a wide Lie subalgebroid of $A$ and $\mathcal{B} := \Gamma_c(B)$. Then $H^G(\mathcal{B})$ is the *minimal integral of* $B$ *over* $G$ in the sense of [24, Thm. 2.3]. This is shown in [31, Prop 3.20].

c) Here is another special case: Let $\mathfrak{g}$ a Lie algebra, $\mathfrak{k}$ a Lie subalgebra, and fix a connected Lie group $G$ integrating $\mathfrak{g}$. Let $\phi\colon K \to G$ be any morphism of Lie groups integrating the inclusion $\iota\colon \mathfrak{k} \hookrightarrow \mathfrak{g}$, where $K$ is assumed to be connected. (For instance, take $K$ to be the connected and simply connected integration of $\mathfrak{k}$.) Then

   $$H(\mathfrak{k}) = K/ker(\phi).$$

   Indeed, since the space of objects of $K$ is just a point, $k_1 \sim k_2$ iff $\phi(k_1) = \phi(k_2)$. Hence $H(\mathfrak{k})$ is a Lie group integrating $\mathfrak{k}$, and the map $\Phi\colon H(\mathfrak{k}) \to G$ induced by $\phi$ is an injective immersion and group homomorphism. In other words, $(H(\mathfrak{k}), \Phi)$ is the Lie subgroup of $G$ integrating $\mathfrak{k}$ (see [29, pp. 92-93]).

We remark that local bisections such as the ones mentioned in the definition of $\mathcal{I}$ in Examples 1.10 a) can be obtained in the following computable way: Put $\mathcal{J}$ the kernel of the $C^\infty(M)$-linear map $d\phi : \Gamma_c(AK) \to \mathcal{B}$. It is an involutive $C^\infty(M)$-submodule of $\Gamma_c(AK)$, whence it corresponds to an involutive $C^\infty(M)$-submodule $\overrightarrow{\mathcal{J}}$ of right-invariant vector fields of $K$. Their associated 1-parameter groups provide bisections which lie in the isotropy of $K$.

## 1.4 Diffeological spaces and their "submersions"

The notion of diffeology was introduced by Souriau in [27] in order to show that the group of diffeomorphisms differentiates to the (infinite-dimensional) Lie algebra of vector fields. To make this article self-contained we give here with an overview of diffeological spaces, following in part the monograph [16].

Let $X$ be a set. We consider the sets of maps:

- $\chi : \mathcal{O}_\chi \to X$, where $\mathcal{O}_\chi$ is an open subset of Euclidean space $\mathbb{R}^k$ for some $k \in \mathbb{N}$. Such a map $\chi$ is called a *$k$-plot*.
- smooth maps $h : \mathcal{O}_h \to \mathcal{O}'_h$, where $\mathcal{O}_h \subseteq \mathbb{R}^n$ and $\mathcal{O}'_h \subseteq \mathbb{R}^k$are open subsets.

We say that a pair $(h, \chi)$ as above is *composable* if $h(\mathcal{O}_h) \subset \mathcal{O}_\chi$. Also, a family $\{\chi_i : \mathcal{O}_i \to X\}_{i\in I}$, where $\mathcal{O}_i \subseteq \mathbb{R}^k$ for every $i \in I$, is called *compatible* if for every $x \in \mathcal{O}_i \cap \mathcal{O}_j$ we have $\chi_i(x) = \chi_j(x)$. In this case there is a unique smallest extension $\chi : \bigcup_{i\in I} \mathcal{O}_i \to X$.

**Definition 1.11.** Consider a set $X$ and for every $k \in \mathbb{N}$ let $\mathcal{P}^k(X)$ be a collection of $k$-plots $\chi$. The collection $\mathcal{P}(X) = \bigcup_{k\in\mathbb{N}} \mathcal{P}^k(X)$ is called a **diffeology** if:

(D1) Every constant map $x : \mathcal{O} \to X$ is in $\mathcal{P}(X)$.

(D2) Given a compatible family $\{\chi_i\}_{i\in I}$ in $\mathcal{P}^k(X)$, its smallest extension is in $\mathcal{P}^k(X)$.

(D3) If $h$ as above and $\chi \in \mathcal{P}^k(X)$ are composable, then $\chi \circ h$ is in $\mathcal{P}^n(X)$.

If $F$ is a set of plots, the intersection of all diffeologies containing $F$ is a diffeology denoted $\langle F\rangle$. It is called the diffeology *generated* by $F$. It is the smallest diffeology containing $F$. In [16, §1.68] it is characterized as follows:

**Proposition 1.12.** Let $F$ be a set of plots. A $k$-plot $\chi : \mathcal{O}_\chi \to X$ belongs to the diffeology $\langle F\rangle$ iff for every point $x$ in $\mathcal{O}_\chi$ there is a neighborhood $\mathcal{O}_x$ of $x$ such that either $\chi|_{\mathcal{O}_x}$ is constant or $\chi|_{\mathcal{O}_x} = \chi' \circ h$, where $\chi' \in F$ and $h$ is such that the pair $(\chi', h)$ is composable.

A diffeology on a set $X$ induces a canonical topology on $X$:

**Definition 1.13.** If $(X, \mathcal{P}(X))$ is a diffeological space, the ***D*-topology** of $X$ is defined as follows: A subset $W$ of $X$ is open iff for each plot $\chi$ in $\mathcal{P}(X)$ the inverse image $\chi^{-1}(W) \subset \mathcal{O}_\chi$ is open.

We also need to specify the notion of a smooth map in this context, and point out some properties.

**Definition 1.14.** Let $(X, \mathcal{P}(X))$ and $(Y, \mathcal{P}(Y))$ two diffeological spaces. A map $f : X \to Y$ is called **smooth** if for every $\chi \in \mathcal{P}(X)$ the composition $f \circ \chi$ is in $\mathcal{P}(Y)$ (It suffices to check this condition for a generating set of $\mathcal{P}(X)$). The collection of smooth maps is denoted $[X, Y]$.

a) Let $(X, \mathcal{P}(X))$ a diffeological space and $X'$ a set. For any map $f : X \to X'$ there exists a finest diffeology $f_*(\mathcal{P}(X))$ on $X'$ which makes $f$ smooth. It is called **pushforward diffeology** and is characterized as follows: A map $\chi' : \mathcal{O}' \to X'$ is a plot in $f_*(\mathcal{P}(X))$ iff for every $r \in \mathcal{O}'$ there exists an open neighborhood $\widetilde{\mathcal{O}'}$ of $r$ such that
   i) either $\chi'|_{\widetilde{\mathcal{O}'}}$ is constant, or
   ii) there exists a plot $\chi : \mathcal{O} \to X$ in $\mathcal{P}(X)$ such that $\chi'|_{\widetilde{\mathcal{O}'}} = f \circ \chi$.

   When $f$ is surjective, this is called **quotient diffeology**, and in the above characterization item i) becomes superflous.

b) In the special case that $X$ is a smooth manifold and $(Y, \mathcal{P}(Y))$ a diffeological space, $f: X \to Y$ is smooth iff for any coordinate neighbourhood $\mathcal{O} \subset X$ we have that $f|_{\mathcal{O}}: \mathcal{O} \to Y$ is a plot in $\mathcal{P}(Y)$. Here and in the following we abuse slightly notation by identifying the coordinate neighbourhood $\mathcal{O}$ with an open subset of $\mathbb{R}^{dim(X)}$ using the chart given by the coordinates. The above is a consequence of the fact that $\mathcal{P}(X)$ is generated by $\{\phi^{-1} : \phi \in \mathfrak{a}\}$ where $\mathfrak{a}$ is an atlas for $M$, see Ex. 1.15 a) below.

The following examples show that the category of diffeological spaces defined above includes manifolds and quotient spaces. Our holonomy groupoid $H^G(\mathcal{B})$ is also an example of a diffeological space, but we will examine it more thoroughly in the next section §2.

**Examples 1.15** (Diffeological spaces and D-topology). a) Let $M$ be an finite dimensional manifold. For every $n \in \mathbb{N}$, we consider $\mathcal{P}^n(M)$ to be the smooth maps $\mathcal{O} \to M$ in the usual sense of differential geometry. If $\mathfrak{a}$ is an atlas of $M$ then the diffeology $\mathcal{P}(M)$ is the one generated by $\mathfrak{a}^{-1} = \{\phi^{-1} : \phi \in \mathfrak{a}\}$. The $D$-topology of $M$ coincides with the topology induced on $M$ by the atlas $\mathfrak{a}$.

If $N$ is another manifold, we have $C^\infty(M, N) = [M, N]$. Notice that the diffeology $\mathcal{P}(M)$ determines uniquely[6] the manifold structure on $M$.

b) Let $\{(X_i, \mathcal{P}(X_i))\}_{i \in I}$ be a family of diffeological spaces, $Y$ a set and $f_i : X_i \to Y$ maps. The *final* diffeology on $Y$ is the one generated by $f_i \circ \chi$ for all $i \in I$ and all plots $\chi \in \mathcal{P}(X_i)$. It is the smallest diffeology making all the $f_i$ smooth.

c) If $\{(X_i, \mathcal{P}(X_i))\}_{i \in I}$ is a family of diffeological spaces, we endow the disjoint union $\coprod_{i \in I} X_i$ with the final diffeology induced by the inclusion maps $\iota_i : X_i \to \coprod_{i \in I} X_i$.

d) Given a diffeological space $(X, \mathcal{P}(X))$ and an equivalence relation $\sim$ on $X$, we endow the quotient $X/\sim$ with the final diffeology induced by the projection $\pi : X \to X/\sim$. It agrees with the quotient diffeology. The $D$-topology of the diffeology on $X/\sim$ coincides with the quotient topology on $X/\sim$ (induced by the projection $\pi$ and the $D$-topology on $X$).

We know from Grothendieck that submersions are a very important class of maps. We argue that in diffeology, this role is played by local subductions defined in [16, §2.17].

**Definition 1.16.** Let $X, Y$ be diffeological spaces. A smooth surjective map $f: X \to Y$ is a **local subduction** if the following holds for every $x \in X$: for any plot $\chi: \mathcal{O}_\chi \to Y$ and any point $r \in \mathcal{O}_\chi$ with $\chi(r) = f(x)$, there is an open neighborhood $V$ of $r$ and a plot $\chi': V \to X$ so that $\chi|_V = f \circ \chi'$ and $\chi'(r) = x$.

**Remarks 1.17.** a) One can regard local subductions as being the analogue – in the diffeological category – of surjective submersions in the category of manifolds. Indeed, when $Y$ is a manifold, the definition of local subduction can be interpreted as saying that every point of $X$ lies on a local section of $f$ (i.e. in the image of a smooth map $\chi'$ such that $f \circ \chi' = Id$).

b) Further, local subductions are open maps with respect to the D-topologies, see [16, §2.20]. Open maps are often regarded as the analogues of submersions in the topological category.

[6] This is due to the fact that, on a smooth manifold $M$, the diffeomorphisms from open subsets of $M$ to open subsets of Euclidean space are exactly the charts of $M$.

Local subductions are in particular **subductions** [16, Ch. 1.46]; the latter are defined as surjective maps $f\colon X \to Y$ such that the diffeology on $Y$ is the quotient diffeology. The property of being a subduction can be phrased as follows: for any plot $\chi$ for $Y$ and any point $r \in \mathcal{O}_\chi$ in its domain, there is an open neighborhood $V$ and a plot $\chi'\colon V \to X$ so that $\chi|_V = f \circ \chi'$.

We finish with a lemma that is not essential for this paper, except for Remark 4.3.

**Lemma 1.18.** Let $(X, \mathcal{P}(X))$ be a diffeological space, and assume the existence of a generating set $\{\chi^i_{op} : \mathcal{O}^i_{op} \to X\}_{i\in I}$ for $\mathcal{P}(X)$ consisting of open maps. Let $\{\chi^j : \mathcal{O}^j \to X\}_{j\in J}$ be any generating set for $\mathcal{P}(X)$. Then for every $j \in J$ there is an open subset $\widetilde{\mathcal{O}^j} \subset \mathcal{O}^j$ such that $\{\chi^j|_{\widetilde{\mathcal{O}^j}}\}_{j\in J}$ is a generating set for $\mathcal{P}(X)$ consisting of maps such that the images $\chi^j(\widetilde{\mathcal{O}^j})$ are open in $X$.

*Proof.* Fix $i \in I$ and $x \in \mathcal{O}^i_{op}$. Since the $\{\chi^j\}_{j\in J}$ form a generating set, there exists a neighborhood $\mathcal{O}^i_{op}(x)$ of $x$, and index $j \in J$ and a smooth map $\tau^x$ making this diagram commute:

$$\begin{array}{ccc} \mathcal{O}^j & \xrightarrow{\chi^j} & X \\ {\scriptstyle \tau^x}\uparrow & \nearrow{\scriptstyle \chi^i_{op}} & \\ \mathcal{O}^i_{op}(x) & & \end{array}$$

Notice that $V^{i,x} := \chi^i_{op}(\mathcal{O}^i_{op}(x))$ is an open subset of $X$. Now for every $j \in J$ define $\widetilde{\mathcal{O}^j} := \cup_{i\in I} \cup_{x\in\mathcal{O}^i_{op}} (\chi^j)^{-1}(V^{i,x})$. One can check that $\{\chi^j|_{\widetilde{\mathcal{O}^j}}\}_{j\in J}$ is a generating set for $\mathcal{P}(X)$. ☐

## 2 Properties of diffeological groupoids

In this section we introduce diffeological groupoids and show that the holonomy groupoid is an example (Prop. 2.3). We also single out in §2.3 some properties that the holonomy groupoid – viewed as a diffeological groupoid – always satisfies: the "holonomy-like" property, the "source-submersive" property, and the "open map/local subduction" property.

Roughly, a diffeological groupoid is a groupoid internal to the category of diffeological spaces, much like a Lie groupoid is essentially a groupoid internal to the category of smooth manifolds. Let us recall the definition from [16, §8.3] in more detail, but imposing the requirement that the set of objects is a manifold:

**Definition 2.1.** A **diffeological groupoid** is a groupoid $H \rightrightarrows M$, where $(H, \mathcal{P}(H))$ is a diffeological space and $M$ a manifold, such that

a) The source and target maps $\mathbf{s}, \mathbf{t} : H \to M$ are smooth.

b) The composition map $H \times_{\mathbf{s},\mathbf{t}} H \to H$ is smooth, when $H \times_{\mathbf{s},\mathbf{t}} H$ is equipped with the subset diffeology of $H \times H$.

c) The inversion map $\iota : H \to H$ is smooth.

d) The unit inclusion map $1 : M \to H$ is smooth.

**Remark 2.2.** Definition 2.1 that diffeological groupoids with a single object ($M = \{pt\}$) are diffeological groups in the sense of [15]. Further, Lie groupoids are special cases of diffeological groupoids. From now on, all diffeological groupoids are assumed to be source-connected. Note that, unlike [10], we do not require the source and target maps to be subductions.

Given two diffeological groupoids $H \rightrightarrows M$ and $H' \rightrightarrows M'$, a **smooth morphism** between them is a groupoid morphism $H \to H'$ which is smooth in the sense of Def. 1.14.

## 2.1 The diffeological structure of the holonomy groupoid

An example of diffeological groupoid is the holonomy groupoid: Let $\mathcal{B}$ be a singular subalgebroid of an integrable Lie algebroid $A = AG$, and $\mathcal{U}^G$ a path-holonomy atlas for $\mathcal{B}$. Recall that $\mathcal{U}^G$ consists of all minimal path-holonomy bisubmersions $(U, \varphi, G)$ together with their inverses and finite compositions. Using $\mathcal{U}^G$ we attach a diffeology to $H^G(\mathcal{B})$ as follows:

- Each $U \in \mathcal{U}^G$ is a smooth manifold, whence a diffeological space with the diffeology described in item (a) of ex. 1.15. Its diffeology is generated by the usual manifold charts $\{(\phi_i^U)^{-1} : \mathcal{O}_i \to U\}$, where $\mathcal{O}_i$ are open subsets of $\mathbb{R}^{dimU}$.
- Put $X = \coprod_{U \in \mathcal{U}^G} U$ and for every $U \in \mathcal{U}^G$ consider the inclusion map $\iota_U : U \to X$. The set $X$ is a diffeological space with the diffeology defined in item (c) of ex. 1.15. Namely its diffeology is generated by $\iota_U \circ (\phi_i^U)^{-1} : \mathcal{O}_i \to X$ for every $U \in \mathcal{U}^G$.
- The *path-holonomy diffeology* is the final diffeology induced by the projection map $\pi : X \to X/\sim$ (see item (d) of ex. 1.15). Its diffeology is generated by $\pi \circ \iota_U \circ (\phi_i^U)^{-1} : \mathcal{O}_i \to H^G(\mathcal{B})$. Notice that $\pi \circ \iota_U$ is the quotient map $q_U : U \to H^G(\mathcal{B})$. Thanks to idem (d) of ex. 1.15, the associated $D$-topology of $H^G(\mathcal{B})$ is the quotient topology.

**Proposition 2.3.** Let $\mathcal{B}$ be a singular subalgebroid of a Lie algebroid $AG$. Then:

a) The path-holonomy diffeology is canonical.

b) The path-holonomy diffeology makes the holonomy groupoid $H^G(\mathcal{B})$ into a diffeological groupoid.

*Proof.* a) To check that the path-holonomy diffeology is canonical, recall that any two path-holonomy atlases are equivalent (see [31, Appendix C], also for the definitions). Let $\mathcal{U}$ and $\mathcal{U}'$ be two equivalent atlases of bisubmersions for the singular subalgebroid $\mathcal{B}$. Denote by $\pi : X \to X/\sim$ and $\pi' : X' \to X'/\sim$ the projection maps, whose final diffeologies on $H^G(\mathcal{B}) = X/\sim = X'/\sim$ we denote by $\mathcal{P}(\mathcal{U})$ and $\mathcal{P}(\mathcal{U}')$ respectively.

If $x \in U \in X$, there is a (smooth) morphism of bisubmersions $f \colon U_0 \to U'$ for some neighborhood $U_0$ of $x$ in $U$ and a bisubmersion $U' \in X'$, since the atlas $\mathcal{U}$ is adapted to $\mathcal{U}'$. Hence $\pi|_{U_0} = \pi' \circ f$ is a smooth map with respect to the diffeology $\mathcal{P}(\mathcal{U}')$. This shows that $\mathcal{P}(\mathcal{U}) \subset \mathcal{P}(\mathcal{U}')$, and by symmetry we obtain the equality.

b) We have to show that the structure maps in Def. 2.1 are smooth. We do so only for b) (the composition map $m$), as the other cases are similar. Let $U_1, U_2 \in \mathcal{U}^G$, and let $\mathcal{O}_1, \mathcal{O}_2$ be coordinate charts there. Then $(q_{U_1}|_{\mathcal{O}_1}, q_{U_2}|_{\mathcal{O}_2})$ is a plot in the diffeology of the Cartesian product $H^G(\mathcal{B}) \times H^G(\mathcal{B})$, and we use the same notation for plot obtained by restriction to the fiber products over $M$. We have to show that $m \circ ((q_{U_1}|_{\mathcal{O}_1}, q_{U_2}|_{\mathcal{O}_2})$ is a plot in the path-holonomy diffeology of $H^G(\mathcal{B})$. This holds because $q_U := m \circ (q_{U_1} \times q_{U_2})$, being the quotient map associated to the composition of two bisubmersions in $\mathcal{U}^G$, lies in $\mathcal{U}^G$.

$$\begin{array}{ccc} U := U_1 \times_{\mathbf{s}_1,\mathbf{t}_2} U_2 & & \\ {\scriptstyle q_{U_1}\times q_{U_2}}\downarrow & \overset{q_U}{\searrow} & \\ H^G(\mathcal{B}) \times_{\mathbf{s},\mathbf{t}} H^G(\mathcal{B}) & \xrightarrow[m]{} & H^G(\mathcal{B}) \end{array}$$

□

**Remarks 2.4.** a) (The path-holonomy diffeology) Let us explain more precisely the diffeological structure of $H^G(\mathcal{B})$ near the identity: Consider a minimal path-holonomy bisubmersion $(U, \varphi, G)$. Shrinking $U$ if necessary, we may use a manifold chart to identify $U$ with an open subset of $\mathbb{R}^{dimU}$. On the other hand, $q_U(U)$ is an open subset of $H^G(\mathcal{B})$ near the identity. Whence the set $H^G(\mathcal{B})_{op} = \bigcup_U q_U(U)$ is an open subset of $H^G(\mathcal{B})$ which contains the identity $M$. (Note that the above union is given only by single path-holonomy bisubmersions, and not by finite compositions of bisubmersions as such.) So the diffeology $\mathcal{Q}$ of $H^G(\mathcal{B})_{op}$ is the one generated by the quotient maps $q_U : U \to H^G(\mathcal{B})_{op}$.

b) (Smoothness of $\Phi$) The groupoid $H^G(\mathcal{B})$ comes together with the map $\Phi : H^G(\mathcal{B}) \to G$. Let us explain the smoothness of $\Phi$ in the diffeological setting: By construction we have that the following diagram commutes, for every minimal path-holonomy bisubmersion $(U, \varphi, G)$:

$$\begin{array}{ccc} U & \xrightarrow{\varphi} & G \\ & {\scriptstyle q_U}\searrow \quad \nearrow{\scriptstyle \Phi} & \\ & H^G(\mathcal{B})_{op} & \end{array} \tag{2.1}$$

This makes the restriction of $\Phi$ to $H^G(\mathcal{B})_{op}$ smooth, *i.e.* an element of $[H^G(\mathcal{B})_{op}, G]$. Allowing finite compositions and inverses in the above discussion, we see that $\Phi \in [H^G(\mathcal{B}), G]$.

c) (The underlying topology) Recall that $H^G(\mathcal{B})$ is a topological groupoid, with the quotient topology induced by $\pi : X \to X/\sim = H^G(\mathcal{B})$. This topology agrees with the $D$-topology underlying the path-holonomy diffeology. This is an immediate consequence of Ex. 1.15 d).

**Remark 2.5.** Given a singular foliation $\mathcal{F}$ that admits a finite resolution, in [19] the authors construct a *holonomy para-Lie* $\infty$*-groupoid* $K_\bullet$, where $\bullet$ ranges over $\mathbb{Z}_{\geqslant 0}$. This object is very close to being a Lie $\infty$-groupoid. They show that the 1-truncation of $K_\bullet$ is precisely the holonomy groupoid $H(\mathcal{F})$ [19, Prop. 2.2]. Indeed, the 1-truncation is the quotient of $K_1$, which consists of bisubmersions, by an equivalence relation induced by the 2-simplices $K_2$. In the proof of [19, Prop. 2.2] the authors show that this equivalence relation coincides with the one of Androulidakis-Skandalis [5], which we recalled at the beginning of §1.3 and which is used to construct the holonomy groupoid. Consequently, the quotient diffeology on the 1-truncation agrees with the path-holonomy diffeology we introduced above (in particular, it makes $H(\mathcal{F})$ into a diffeological groupoid).

However, as it is pointed out in [19, Rem. 1.13], not all smooth singular foliations in the sense of [5] admit a finite resolution. On the other hand, the results in this sequel apply to *any* smooth singular foliation, including the (non-analytic) cases mentioned in [19, Rem. 1.13] which do not admit a finite resolution.

## 2.2 Bisections of diffeological groupoids

Let $H \rightrightarrows M$ be a diffeological groupoid. The key to the differentiation we carry out in §3 is the notion of a bisection, which we discuss here.

**Definition 2.6.** a) A **global bisection** of $H$ is a smooth map $\mathbf{b} : M \to H$ which is right-inverse to $\mathbf{s}_H$ and such that $\mathbf{t}_H \circ \mathbf{b}$ is a diffeomorphism onto its image. We use the term **bisection** when the domain is any open subset of $M$.

b) Let $I$ be an interval. Consider a family $\{\mathbf{b}_\lambda : V \to H\}_{\lambda \in I}$ of bisections. This family is called **smooth** iff the map $V \times I \to H, (x, \lambda) \mapsto \mathbf{b}_\lambda(x)$ is smooth as a map between diffeological spaces.

**Remarks 2.7.** a) Fix a generating set of plots $\mathcal{X} = \{\chi : \mathcal{O}_\chi \to H\}$ in $\mathcal{P}(H)$. If $\mathbf{b} : V_{\mathbf{b}} \to H$ is a bisection and $V_{\mathbf{b}}$ is small enough, there exists a plot $\chi : \mathcal{O}_\chi \to H$ in $\mathcal{X}$ and a smooth map $\mathbf{u} : V_{\mathbf{b}} \to \mathcal{O}_\chi$ satisfying $\chi \circ \mathbf{u} = \mathbf{b}$. This follows from item c) in §1.4 and Prop. 1.12.

b) An example of a global bisection is the identity bisection $1_M$ of $H$, namely the unit inclusion map, which is smooth by the definition of diffeological groupoid. For the holonomy groupoid $H^G(\mathcal{B})$, this fact can be also seen as follows: Let $(U, \mathbf{t}, \mathbf{s})$ be a minimal path-holonomy bisubmersion and put $V = \mathbf{s}(U)$. Since $V \to U, x \mapsto (x, 0)$ is a bisection of the bisubmersion, it follows that $1_V : V \to H^G(\mathcal{B})$ defined by $1_V(x) = q(x, 0)$ is a bisection of the holonomy groupoid.

Global bisections can be multiplied: for $i = 1, 2$ consider global bisections $\mathbf{b}_i : M \to H$. We define $\mathbf{b}_2 * \mathbf{b}_1 : M \to H$ by

$$(\mathbf{b}_2 * \mathbf{b}_1)(x) = \mathbf{b}_2((\mathbf{t}_H \circ \mathbf{b}_1)(x)) \cdot \mathbf{b}_1(x)$$

where the product on the right-hand side is the groupoid multiplication of $H$. Similarly, the inverse to a global bisection $\mathbf{b}$ is the global bisection $\mathbf{b}^{-1}$ determined by $\mathbf{b}^{-1}((\mathbf{t} \circ \mathbf{b})(x)) = (\mathbf{b}(x))^{-1}$ for all $x \in M$ (see [23, Prop. 1.42]).

**Definition 2.8.** Let $I$ be an open interval containing 0. A smooth family of global bisections $\{\mathbf{b}_t\}_{t \in I}$ of $H$ is called a **1-parameter group** iff

a) $\mathbf{b}_0 = 1_M$,

b) $\mathbf{b}_{\lambda+\mu} = \mathbf{b}_\lambda * \mathbf{b}_\mu$ whenever $\lambda, \mu, \lambda + \mu \in I$.

The next result proves the existence of many 1-parameter groups of global bisections for $H^G(\mathcal{B})$.

**Proposition 2.9.** Let $\mathcal{B}$ be a singular subalgebroid. Then for every $x \in M$ there is a neighborhood $V$ with the following property: every $\boldsymbol{\alpha} \in \mathcal{B}$ with support in $V$ gives rise to an 1-parameter group of compactly supported[7] global bisections $\{\mathbf{b}_\lambda\}$ of $H^G(\mathcal{B})$, satisfying $\frac{d}{d\lambda}|_{\lambda=0}(\Phi \circ \mathbf{b}_\lambda) = \boldsymbol{\alpha}$

*Proof.* Given $x \in M$, consider a minimal path-holonomy bisubmersion $(U, \varphi, G)$ at $x$, and define $V := (\mathbf{s} \circ \varphi)(U)$, which is open in $M$ by Prop. 1.9. Pick a section $\boldsymbol{\alpha} \in \mathcal{B}$ with support in $V$, and consider the corresponding right-invariant vector field $\overrightarrow{\boldsymbol{\alpha}}$ of $G$. From definition 1.3 we have $\varphi^{-1}(\overrightarrow{\mathcal{B}}) = \Gamma_c(U; \ker d\mathbf{s}_U)$, whence there is a vector field $\xi \in \Gamma(U; \ker d\mathbf{s}_U)$ which is $\varphi$-related to $\overrightarrow{\boldsymbol{\alpha}}$.

[7] More precisely this means: the bisections agree with the identity bisection of $V$ outside a compact subset of $V$.

Looking at [23, §3.6] we may assume that there is an interval $I = (-\epsilon, \epsilon)$ such that the flow $\{\mathbf{A}_\lambda\}_{\lambda\in I}$ of $\vec{\boldsymbol{\alpha}}$ is defined on an open subset $W$ of $G$ which contains an open neighborhood of $V$ in the base manifold $M$. Since $\xi$ is $\varphi$-related to $\vec{\boldsymbol{\alpha}}$ we can assume that the flow $\{\boldsymbol{\Xi}_\lambda\}_{\lambda\in I}$ of $\xi$ is defined in an open subset $\hat{W}$ of $U$ which contains $V \times \{0\}$. The situation is summarized in the following diagram:

$$\begin{array}{ccc} U & \xrightarrow{\varphi} & G \\ \xi, \{\boldsymbol{\Xi}_\lambda\} & & \vec{\boldsymbol{\alpha}}, \{\mathbf{A}_\lambda\} \end{array}$$

For every $\lambda$ define $\mathbf{u}_\lambda : V \to U$ by

$$\mathbf{u}_\lambda(x) = \Xi_\lambda(x, 0). \tag{2.2}$$

The family $\{\mathbf{u}_\lambda\}_{\lambda\in I}$ is obviously smooth. Composing each $\mathbf{u}_\lambda$ with the quotient map $q_U : U \to H^G(\mathcal{B})$ we obtain a 1-parameter family of bisections $\{\mathbf{b}_\lambda\}_{\lambda\in I}$ of $H^G(\mathcal{B})$, with the following properties:

- It satisfies $\frac{d}{d\lambda}|_{\lambda=0}(\Phi \circ \mathbf{b}_\lambda) = \boldsymbol{\alpha}$ by construction, using that $\Phi \circ q_U = \varphi$.
- The bisections $\mathbf{b}_\lambda$ are compactly supported. Indeed, for all $y \in V$ with $y \notin Supp(\boldsymbol{\alpha})$ we have $\mathbf{b}_\lambda(y) = 1_y \in H^G(\mathcal{B})$, for every $\lambda$. This follows from the second bullet point in §1.3 and the fact that in a neighborhood of $y$, $\mathbf{u}_\lambda$ is a bisection of $U$ which carries the identity bisection of $G$ (i.e. $\varphi(\mathbf{u}_\lambda) \subset 1_M$ in such a neighborhood).
- The family $\{\mathbf{b}_\lambda\}_{\lambda\in I}$ is a 1-parameter group, since it satisfies the properties of definition 2.8:
  a) $\mathbf{b}_0 \subset 1_M$, as a consequence of $\mathbf{u}_0(x) = (x, 0)$.
  b) We have to prove that the bisections $\mathbf{b}_\lambda * \mathbf{b}_\mu$ and $\mathbf{b}_{\lambda+\mu}$ of $H^G(\mathcal{B})$ coincide. They both carry the bisection
  $$\mathbf{A}_\lambda * \mathbf{A}_\mu = \mathbf{A}_{\lambda+\mu}$$
  of $G$, since $\Phi$ is a groupoid morphism, since $\Phi \circ q_U = \varphi \colon U \to G$ and $\xi$ is $\phi$-related to $\vec{\boldsymbol{\alpha}}$. Here the equality holds because $\{\mathbf{A}_\lambda\}$ is the flow of a right invariant vector field on $G$ (namely, $\vec{\boldsymbol{\alpha}}$). Since both $\mathbf{b}_\lambda * \mathbf{b}_\mu$ and $\mathbf{b}_{\lambda+\mu}$ carry the same bisection of $G$, they coincide, by the second bullet point in §1.3.

□

## 2.3 Special properties of diffeological groupoids

For diffeological groupoids, we introduce three properties, which are always satisfied by the holonomy groupoid and by Lie groupoids.

1) The *holonomy-like* property, which plays an important role in the differentiation of diffeological groupoids explained in §3.1.
2) The *source-submersive* property, which implies that the source map and target map are both local subductions, thus open maps.
3) A third property which comes in two flavours, the *open map* and the *local subduction* property, which will be important for the notion of integral in §4.

Property 1) seems to be independent of the other two (see for instance Lemma 5.4 together with Prop. 5.6).

Property 3) is stronger than property 2): indeed the following chain of implications holds, by Prop. 2.17 and Remark 2.23:

$$\text{local subduction property} \;\Rightarrow\; \text{open map property} \;\Rightarrow\; \text{source-submersive property.}$$

**Remark 2.10** (A smoothness property)**.** Recall that the definition a Lie groupoid does not only entail that the spaces and maps involved are smooth, but also that the source and target maps are submersions. We interpret the third property above as a "smoothness condition" on diffeological groupoids, analog to the submersion condition for the source and target of Lie groupoids. An example of diffeological groupoid that does not satisfies the third condition is the graph of a singular subalgebroid, see Prop. 5.6.

Before introducing the three properties, we present an important lemma, which applies to all diffeological groupoids. It will be essential in the proofs of Prop. 2.17 and Lemma A.2.

**Lemma 2.11.** Let $H \rightrightarrows M$ be a diffeological groupoid. Let $\{\chi : \mathcal{O}_\chi \to H\}$ be any generating set of plots in $\mathcal{P}(H)$. Then for every $x \in M$ there is a neighbourhood $V \subset M$, a plot $\chi$ in the generating set, and a submanifold $e \subset \mathcal{O}_\chi$, such that

i) $\chi|_e \colon e \to 1(V)$ is a diffeomorphism,

ii) $\mathbf{s}_H \circ \chi$ and $\mathbf{t}_H \circ \chi$ are submersions at points of $e$.

*Proof.* Fix $x \in M$. The inclusion of the identity elements $1 : M \to H$ is smooth by Def. 2.1. There is an open neighbourhood $V$ of $x$ in $M$ so that $1|_V$ is a plot for $H$ in $\mathcal{P}(H)$, since $M$ is a manifold, see item c) in §1.4 (after identifying $V$ with an open subset of $\mathbb{R}^{dim(M)}$ by means of a chart). There exists a plot $\chi : \mathcal{O}_\chi \to H$ in the generating set and a smooth map $\tilde{1} \colon V \to \mathcal{O}_\chi$ with $\chi \circ \tilde{1} = 1|_V$, by Prop. 1.12, after shrinking $V$ if necessary.

Notice that, since $\mathbf{s}_H \circ 1|_V = Id_V$, the map $\tilde{1}$ is a smooth section of $\mathbf{s}_H \circ \chi \colon \mathcal{O}_\chi \to V$. This has two consequences.

The first is that $e := \tilde{1}(V) \subset \mathcal{O}_\chi$ is a submanifold of $\mathcal{O}_\chi$, shrinking $V$ if necessary, because $\tilde{1}$ is an injective immersion. Further $\chi(e) = 1(V)$ is an open neighbourhood of $1_x$ in $1_M$. Also, $\chi|_e \colon e \to \chi(e) = 1(V)$ is a diffeomorphism, because it is the inverse of the diffeomorphism $\tilde{1} \colon V \to e$.

$$\begin{array}{ccc} e \subset \mathcal{O}_\chi & \xrightarrow{\;\chi\;} & H \\ & \nwarrow^{\tilde{1}} & \uparrow{\scriptstyle 1|_V} \\ & & V \end{array} \tag{2.3}$$

The second consequence is item ii) in the statement. ☐

### 2.3.1 The "holonomy-like" property

This property of diffeological groupoids will be used in §3.1. It is a reasonable property in view of Lemma 2.11, and loosely speaking it implies that there is a generating set for $H$ nearby the identity section consisting of "few plots" (see Prop. 1.12). The plots of this generating set can be heuristically described as being "semi-local subductions over points of $1_M$", a variation of the notion of local subduction where a condition is imposed on submanifolds rather than points.

**Definition 2.12.** A diffeological groupoid $H \rightrightarrows M$ is called **holonomy-like** if there exists an open neighborhood $\overset{\circ}{H}$ of $M$ in $H$ and a generating set of plots $\mathcal{X} = \{\chi : \mathcal{O}_\chi \to \overset{\circ}{H}\}$ in $\mathcal{P}(\overset{\circ}{H})$ such that the following holds for all $\chi, \chi' \in \mathcal{X}$:

- Let $e \subset \mathcal{O}_\chi, e' \subset \mathcal{O}_{\chi'}$ be submanifolds such that $\chi(e) = \chi'(e')$ is an open subset of $1_M \subset H$ and $\chi|_e : e \to \chi(e)$ and $\chi'|_{e'} : e' \to \chi'(e')$ are diffeomorphisms. Let $x' \in e'$. Then there exists a smooth map $k : \mathcal{O}_{\chi'} \to \mathcal{O}_\chi$ with $k(e') = e$ and $\chi \circ k = \chi'$, after shrinking $\mathcal{O}_{\chi'}$ to a smaller neighborhood of $x'$ if necessary.

$$\begin{array}{ccc} & H & \\ {\scriptstyle \chi'} \nearrow & & \nwarrow {\scriptstyle \chi} \\ \mathcal{O}_{\chi'} & \overset{k}{\dashrightarrow} & \mathcal{O}_\chi \end{array} \tag{2.4}$$

**Proposition 2.13.** The holonomy groupoid $H^G(\mathcal{B})$ is holonomy-like.

*Proof.* Take $\overset{\circ}{H} = H^G(\mathcal{B})$ and $\mathcal{X}$ equal to the path holonomy atlas (see §1.2). Let $\chi, \chi', e, e'$ be as in Def. 2.12. In particular, $e$ is a bisection of the bisubmersion $\mathcal{O}_\chi$, and similarly for $e'$ and $\mathcal{O}_{\chi'}$. It follows that $(\Phi \circ \chi)(e)$ and $(\Phi \circ \chi')(e')$ are both mapped into the identity bisection of $G$. Hence there is a morphism of bisubmersions $k : \mathcal{O}_{\chi'} \to \mathcal{O}_\chi$ with $k(e') = e$, as follows from the proof of [31, Prop. 3.2, Cor. 3.3] (that proof is analogous to [5, Prop. 2.10 b), Cor. 2.11 a]). We have $\chi \circ k = \chi'$ by the definition of the equivalence relation used in the construction of the holonomy groupoid (see §1.3), since $k$ is a morphism of bisubmersions. □

**Remark 2.14.** Lie groupoids are holonomy-like. This is immediate since the diffeology of a Lie groupoid is generated by inverses of charts (see Ex. 1.15 a)).

### 2.3.2 The "source-submersive" property

Recall that the source map of a Lie groupoid is a surjective submersion. An analogous statement holds for diffeological groupoids satisfying the following property:

**Definition 2.15.** A diffeological groupoid $H \rightrightarrows M$ is **source-submersive** if through every point of $H$ there is a (local) bisection.

Indeed, in this case the maps $\mathbf{s} : H \to M$ and $\mathbf{t} : H \to M$ are local subductions, and in particular open maps (see Remark 1.17).

The assumptions of the next lemma, and also those of Prop. 2.17, are satisfied by the holonomy groupoid $H^G(\mathcal{B})$ (thanks to Prop. 1.9) and by Lie groupoids, which therefore are source-submersive.

**Lemma 2.16.** Let $H \rightrightarrows M$ be a (source-connected) diffeological groupoid. Assume there exist plots $\{\chi_i : \mathcal{O}_{\chi_i} \to H\}$ in $\mathcal{P}(H)$ such that

i) $\mathbf{s} \circ \chi_i : \mathcal{O}_{\chi_i} \to M$ and $\mathbf{t} \circ \chi_i$ are submersions, for all $i$,

ii) $\overset{\circ}{H} := \cup_i Image(\chi_i)$ is an open neighborhood of $M$ in $H$.

Then $H$ is source-submersive.

*Proof.* Claim: *Through every point of $\overset{\circ}{H}$ there is a bisection.*

Let $h \in \overset{\circ}{H}$. Then $h = \chi_i(x)$ for some $x \in \mathcal{O}_{\chi_i}$. By assumption i), there is a smooth submanifold of $\mathcal{O}_{\chi_i}$ through $x$ which is transverse to both the $\mathbf{s} \circ \chi_i$-fibers and the $\mathbf{t} \circ \chi_i$-fibers. It is the image of a unique local section $\sigma$ of $\mathbf{s} \circ \chi_i$. Hence $\chi_i \circ \sigma$ is a smooth bisection of $H$ through $h$, proving the claim.

$$\begin{array}{ccc} \mathcal{O}_{\chi_i} & \xrightarrow{\chi_i} & H \\ & {\nwarrow}^{\sigma} & \downarrow \mathbf{s} \\ & & M \end{array} \tag{2.5}$$

△

Every source-connected topological groupoid is generated by any symmetric open neighborhood of the identity section (this can be proven as in [23, Prop. 1.5.8]). Hence, by assumption ii), any point $h \in H$ is the product $h_1 \cdots h_k$ of points in $\overset{\circ}{H}$. By the claim, there is a smooth bisection $\mathbf{b}_j$ through $h_j$ for all $j$. Their product $\mathbf{b}_1 * \cdots * \mathbf{b}_k$ is a smooth bisection through $h$. This proves the statement. □

The following proposition provides a sufficient criterion for being source-submersive, and an additional property. The notion of diffeology generated by open maps appearing here will be formalized in §2.3.3.

**Proposition 2.17.** Let $H \rightrightarrows M$ be a diffeological groupoid. Assume there is an open neighborhood $\overset{\circ}{H}$ of $1_M$ in $H$ whose diffeology is generated by *open maps*, i.e.: there exists a generating set of plots $\{\chi_i : \mathcal{O}_{\chi_i} \to \overset{\circ}{H}\}$ of $\mathcal{P}(\overset{\circ}{H})$ such that each $\chi_i$ is an open map. Then:

i) the diffeological groupoid $H$ is source-submersive,

ii) there is a neighborhood $\breve{H}$ of the identity $1_M$ in $H$ such that for all $h \in \breve{H}$ there is a 1-parameter family of global bisections $\{\mathbf{b}_\lambda\}_{\lambda \in [0,1]}$ with $\mathbf{b}_0 = 1_M$ and $h \in \mathbf{b}_1$.

*Proof.* For every $x \in 1_M$, by Lemma 2.11 we know that there is a neighbourhood $V$ of $x$ in $M$, a plot $\chi \colon \mathcal{O}_\chi \to \overset{\circ}{H}$ in the above generating set, and a submanifold $e$ of $\mathcal{O}_\chi$, such that

- $\chi|_e \colon e \to 1(V) \subset 1_M$ is a diffeomorphism,
- $\mathbf{s}_H \circ \chi$ and $\mathbf{t}_H \circ \chi$ are a submersion at points of $e$.

Therefore there is a neighborhood $\overset{\circ}{\mathcal{O}}_\chi$ of $e$ in $\mathcal{O}_\chi$ so that

- $\mathbf{s}_H \circ \chi$ and $\mathbf{t}_H \circ \chi$ are submersions when restricted to $\overset{\circ}{\mathcal{O}}_\chi$,
- the image of $\overset{\circ}{\mathcal{O}}_\chi$ under $\chi$ is a neighborhood of $1(V)$ in $H$ (here we use the assumption that $\chi$ is an open map).

We now prove the two items of the lemma.

i) The set of plots $\{\chi_i|_{\overset{\circ}{\mathcal{O}}_{\chi_i}} : \overset{\circ}{\mathcal{O}}_{\chi_i} \to \overset{\circ}{H}\}$ satisfies the assumptions of Lemma 2.16.

ii) <u>Claim:</u> *Fix an index* $i$*. Shrinking* $\mathring{\mathcal{O}}_{\chi_i}$ *to a smaller neighborhood of* $e_i$ *if necessary, every point of* $p \in \mathring{\mathcal{O}}_{\chi_i}$ *lies on the image of* $\mathbf{u}_1$ *for some family* $\{\mathbf{u}_\lambda\}$ *of sections of* $\mathbf{s}_H \circ \chi_i \colon \mathring{\mathcal{O}}_{\chi_i} \to V_i$ *depending smoothly on* $\lambda \in [0,1]$*, transverse to the fibers of* $\mathbf{t}_H \circ \chi_i$*, with compact support, and with* $\mathbf{u}_0 = (\chi_i|_{e_i})^{-1}$.

Given the claim, we have that $\chi_i(p)$ lies on the image of $\mathbf{b}_1$, where

$$\mathbf{b}_\lambda := \chi_i \circ \mathbf{u}_\lambda, \qquad \lambda \in [0,1],$$

is a 1-parameter family of bisections in $H$ with support in $V_i$, which hence can be extended trivially to global bisections. Taking $\breve{H} = \cup_i \chi_i(\mathring{\mathcal{O}}_{\chi_i})$, this proves the desired statement.

To show the claim, since we are working locally, we may assume that $\mathbf{s}_H \circ \chi_i \colon \mathring{\mathcal{O}}_{\chi_i} \to V_i$ is a neighborhood of the zero section of a metric vector bundle with zero section $e_i$. Clearly $e_i$ is transverse to the fibers of $\mathbf{t}_H \circ \chi_i$. For every unit length vector $v \in \mathring{\mathcal{O}}_{\chi_i}$, extending $v$ to a compactly supported section of this vector bundle, we obtain a vertical fiberwise constant vector field, and there is $\epsilon > 0$ so that the image of $e_i$ under its flow is transverse to the fibers of $\mathbf{t}_H \circ \chi_i$ for all times in $[0,\epsilon)$. Hence each of the points $tv$, for $t \in [0,\epsilon)$, satisfies the statement of the claim.

□

### 2.3.3 The "open map" and "local subduction" property

Here we look at diffeological spaces and groupoids admitting a generating set of plots which define either an open map, or a local subduction. Recall these two kinds of maps are the analogues of (surjective) submersions respectively in the topological and the diffeological category.

**Definition 2.18.** Let $(X, \mathcal{P}(X))$ a diffeological space. We say that the diffeology is **generated by open maps** if there exists a generating set $\{\chi_i \colon \mathcal{O}_i \to X\}_{i \in I}$ of plots in $\mathcal{P}(X)$ such that $\chi_i \colon \mathcal{O}_i \to X$ is an open map for all $i$, w.r.t. the D-topology. Equivalently, so that the combined map $\sqcup \chi_i \colon \sqcup \mathcal{O}_i \to X$ is an open map.

**Example 2.19.** An example of diffeological space whose diffeology is *not* generated by open maps is the following. Let $\gamma \colon \mathbb{R} \to \mathbb{R}^2$ be a smooth map which is injective except at the points 0 and 1, for which we have $\gamma(0) = \gamma(1)$. Let $X = \gamma(\mathbb{R})$ be the image of this map (a curve in the plane with one self-intersection), endowed with the quotient diffeology. Since $Id_\mathbb{R}$ is a generating plot for $\mathbb{R}$, the map $\gamma \colon \mathbb{R} \to X$ is a generating plot for $X$. We know that the D-topology of $X$ is the quotient topology induced by the map $\gamma \colon \mathbb{R} \to X$, see Examples 1.15. The plot $\gamma$ is not an open map: taking $0 < \epsilon < 1$, we have that $\gamma((-\epsilon,\epsilon))$ is not an open subset of $X$ in the quotient topology, because its pre-image under $\gamma$ is $(-\epsilon,\epsilon) \sqcup \{1\}$, which is not open in $\mathbb{R}$.

We already applied this notion to a diffeological groupoid $H$ in Prop. 2.17, where we showed that if the diffeology of $H$ is generated by open maps, then $H$ is source-submersive.

**Lemma 2.20.** Let $H$ be a diffeological groupoid. The following are equivalent:

- the diffeology of $H$ is generated by open maps,

- there is an open neighborhood $\overset{\circ}{H}$ of $1_M$ whose diffeology is generated by open maps.

*Proof.* The implication from top to bottom is easy. For the other implication, notice that the existence of such $\overset{\circ}{H}$ implies that a generating set for $\mathcal{P}(H)$ is obtained composing the plots $\{\chi_i : \mathcal{O}_{\chi_i} \to \overset{\circ}{H}\}$ as in Prop. 2.17 with the right-translations $R_{\mathbf{b}}$, where the $\mathbf{b}$ range over all bisections of $H$. Each $R_{\mathbf{b}} \circ \chi_i$ is an open map, being a composition of such. □

**Definition 2.21.** Let $(X, \mathcal{P}(X))$ a diffeological space. We say that the diffeology is **generated by local subductions** if there exists a generating set $\{\chi_i \colon \mathcal{O}_i \to X\}_{i \in I}$ of plots in $\mathcal{P}(X)$ such that $\sqcup\chi_i \colon \sqcup \mathcal{O}_i \to X$ is a local subduction.

Recall that local subductions were defined in Def. 1.16, are surjective by definition, and are always open maps. In the following remark we rephrase Def. 2.21 in terms of the individual plots $\chi_i$.

**Lemma 2.22.** Given a diffeological space $X$ and a generating set of plots $\{\chi_i \colon \mathcal{O}_i \to X\}_{i \in I}$ in $\mathcal{P}(X)$, the following is equivalent:

- the map $\sqcup\chi_i$ is a local subduction,
- the following three conditions are satisfied:
  - i) for each $i \in I$, the image $\chi_i(\mathcal{O}_i)$ is an open subset of $X$ (where $X$ is endowed with the D-topology),
  - ii) for each $i \in I$, the map $\chi_i \colon \mathcal{O}_i \to \chi_i(\mathcal{O}_i)$ is a local subduction, where $\chi_i(\mathcal{O}_i)$ is endowed with the subspace diffeology,
  - iii) the map $\sqcup\chi_i \colon \sqcup \mathcal{O}_i \to X$ is surjective.

*Proof.* We argue as follows. Assuming the first condition, the second is clear. For the other implication, let $\chi \colon \mathcal{O} \to X$ be any plot in $\mathcal{P}(X)$, $r \in \mathcal{O}$, and $y \in \mathcal{O}_i$ (for some $i \in I$) such that $\chi(r) = \chi_i(y)$. Since $\chi_i(\mathcal{O}_i)$ is an open subset of $X$, there is a neighborhood $V \subset \mathcal{O}$ of $r$ such that $\chi|_V$ takes values in $\chi_i(\mathcal{O}_i)$. The fact that $\chi_i \colon \mathcal{O}_i \to \chi_i(\mathcal{O}_i)$ is a local subduction means that, shrinking $V$ if necessary, there is a smooth map $h \colon V \to \mathcal{O}_i$ such that $\chi_i \circ h = \chi|_V$ and $h(r) = y$. Viewing $h$ as a map $V \to \sqcup\mathcal{O}_i$ finishes the argument. □

Notice that conditions ii) and iii) above do not imply condition i): take the diffeology on $X = \mathbb{R}$ generated by the plots $\chi_1 = x^2 \colon \mathbb{R} \to \mathbb{R}$ and $\chi_2 = -x^2 \colon \mathbb{R} \to \mathbb{R}$. Each $\chi_i$ is a local subduction, but $\chi_1(\mathbb{R}) = \mathbb{R}_{\geqslant 0}$ is not open in the D-topology, because its preimage under $\chi_2$ is the singleton $\{0\}$.

The relation between Def. 2.18 and Def. 2.21 is manifest, because local subductions are always open maps:

**Remark 2.23.** Let $(X, \mathcal{P}(X))$ a diffeological space. If the diffeology is generated by local subductions, then it is generated by open maps.

The following proposition states that the quotient map from a path-holonomy atlas to the holonomy groupoid plays the role of a surjective submersion in the diffeological category.

**Proposition 2.24.** Let $\mathcal{B}$ be a singular subalgebroid of a Lie algebroid $AG$, and let $\mathcal{U}^G$ be a path-holonomy atlas for $\mathcal{B}$. The quotient map $\natural : \coprod_{U \in \mathcal{U}^G} U \to H^G(\mathcal{B})$ is a local subduction.

*Proof.* The map $\natural$ is surjective, and it is smooth because the path-holonomy diffeology is exactly the quotient diffeology induced by $\natural$. Fix a point $u$ in a bisubmersion $U \in \mathcal{U}^G$, denote $h := \natural(u)$ the element of the holonomy groupoid it represents. Take any plot $\chi\colon \mathcal{O}_\chi \to H^G(\mathcal{B})$ and any point $r \in \mathcal{O}_\chi$ with $\chi(r) = h$. There is a connected open neighborhood $V$ of $r$ and a plot $\chi'\colon V \to \coprod_{U \in \mathcal{U}^G} U$ so that $\chi|_V = \natural \circ \chi'$, by definition of quotient diffeology. The image $u' := \chi'(r)$ lies in some bisubmersion $U'$. The point $u'$ might differ from the point $u$ we fixed, but we know that $\natural(u') = h = \natural(u)$. By the definition of holonomy groupoid, this means that there is a morphism of bisubmersions $\phi\colon U' \to U$ mapping $u'$ to $u$, shrinking $U'$ if necessary. As for every morphism of bisubmersions, we have $\natural \circ \phi = \natural$. Thus the composition $\phi \circ \chi'\colon V \to U$ (which is well-defined after shrinking $V$ if necessary) is a plot mapping $r$ to $u$ and lifting the plot $\chi|_V$.

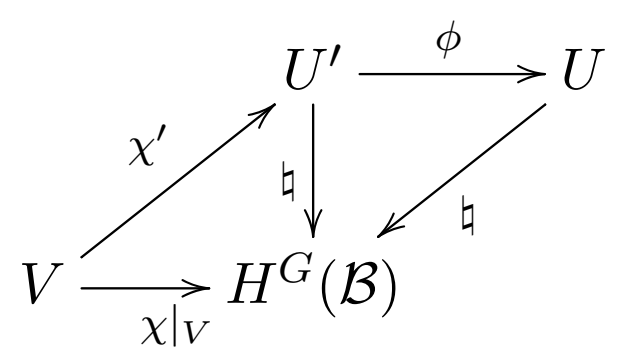

□

Prop. 2.24 immediately implies:

**Corollary 2.25.** The diffeology of the holonomy groupoid is generated by local subductions.

**Remark 2.26.** The diffeology of any manifold, thus also that of any Lie groupoid, is generated by local subductions. This holds because the diffeology is generated by the inverses of manifold charts (see Examples 1.15 a)), and using Lemma 2.22.

# 3 Global differentiation

We explain how a diffeological groupoid together with a morphism to a Lie groupoid $G$ – both satisfying certain properties – gives rise to a singular subalgebroid (Thm. 3.2). We formalize this in the definition of differentiation (Def. 3.3). The primary example is the following: the holonomy groupoid of a singular subalgebroid $\mathcal{B}$ differentiates to $\mathcal{B}$ itself (Thm. A.7). The proofs are all collected in Appendix A, where we also address the functoriality of differentiation.

## 3.1 Differentiation of global bisections to singular subalgebroids

In this subsection, in view of the holonomy groupoid, we fix the following data:

- a diffeological groupoid $H \rightrightarrows M$, which is holonomy-like (see Def. 2.12),
- a Lie groupoid $G$ and a smooth morphism of diffeological groupoids $\Psi : H \to G$ covering $Id_M$.

We use the apparatus developed in §2.1 and §2.2 to specify a differentiation process for families of global bisections which gives rise to a singular subalgebroid of $AG$ (Thm. 3.2 and Def. 3.3 below).

To this end, first note that if $\{\mathbf{b}_\lambda\}_{\lambda\in I}$ is a smooth family (resp. 1-parameter group) of global bisections of $H$ then $\{\Psi\circ\mathbf{b}_\lambda\}_{\lambda\in I}$ is a smooth family (resp. 1-parameter group) of global bisections of $G$. Now let us consider the following set of "velocity vectors":

$$\mathcal{S} := \left\{ \frac{d}{d\lambda}|_{\lambda=0}(\Psi\circ\mathbf{b}_\lambda)\colon\ \{\mathbf{b}_\lambda\}_{\lambda\in I} \text{ family of global bisections for } H \text{ s.t. } \mathbf{b}_0 = Id_M \right\} \cap \Gamma_c(AG).$$

We will argue that $\mathcal{S}$ "contains many elements" in Prop. A.5 below, by constructing elements out of 1-parameter families of local bisections.

Let us make the following assumption, which roughly says that the elements of $\mathcal{S}$ arise from 1-parameter *groups*. In §A.3 we show that the holonomy groupoid satisfies it.

**Assumption 3.1.** For every $x\in M$ there is a neighborhood $V$ with this property:

$$\{\boldsymbol{\alpha}\in\mathcal{S} : Supp(\boldsymbol{\alpha})\subset V\} \quad \subset$$
$$\left\{ \frac{d}{d\lambda}|_{\lambda=0}(\Psi\circ\mathbf{b}_\lambda)\colon\ \{\mathbf{b}_\lambda\}_{\lambda\in I} \text{ 1-parameter } \textit{group} \text{ of global bisections for } H \right\}.$$

The main statement of this section is:

**Theorem 3.2.** Let $H\rightrightarrows M$ and $\Psi : H\to G$ be as above (in particular $H$ is holonomy-like), satisfying Assumption 3.1. Then the set $\mathcal{S}\subseteq\Gamma_c(AG)$ is a singular subalgebroid of $AG$.

We refer the reader to Appendix A for the proof of Theorem 3.2, examples, and observations on the fuctoriality of differentiation.

Now the next definition seems in order, in view of Theorem 3.2:

**Definition 3.3 (Differentiation).** Let $H\rightrightarrows M$ be a diffeological groupoid, $G\rightrightarrows M$ a Lie groupoid and $\Psi : H\to G$ a smooth morphism of diffeological groupoids covering $Id_M$. If

- $H$ is holonomy-like (Def. 2.12), and
- this data satisfies Assumption 3.1,

then we say that $(H,\Psi)$ **differentiates** to the singular subalgebroid $\mathcal{S}$ of $\Gamma_c(AG)$.

For holonomy groupoids and Lie groupoids the differentiation procedure yields the following, as we prove in Thm. A.7 and Prop. A.6:

- Given a singular subalgebroid $\mathcal{B}$, its holonomy groupoid $H^G(\mathcal{B})$ differentiates to $\mathcal{B}$
- If $H$ is a source-connected Lie groupoid and $\Psi\colon H\to G$ a Lie groupoid morphism covering $Id_M$, then $H$ differentiates to $\Psi_*(\Gamma_c(AH))$.

We finish with a few comments on the hypotheses of Theorem 3.2.

**Remarks 3.4.** a) In Theorem 3.2, the property of $H$ being holonomy-like (*cf.* definition 2.12) is used only to prove that the module $\mathcal{S}$ is locally finitely generated, in Lemma A.2.

There do exist interesting examples of involutive $C^\infty$-submodules of vector fields which are not finitely generated. (For instance, consider the partition to $\mathbb{R}$ to $(0, +\infty)$, $\{0\}$ and $\{x\}$ for every $x < 0$ and take the $C_c^\infty(\mathbb{R})$-module generated by all vector fields on $\mathbb{R}$ which are tangent to these submanifolds.) We do not know if diffeological groupoids can be used to treat the differentiation/integration process of non locally finitely generated modules of vector fields.

b) Assumption 3.1 is used only to prove that $\mathcal{S}$ is involutive, in Lemma A.3. The interested reader may check that a variation of the proof of Lemma A.3 gives the following result; instead of making Assumption 3.1, suppose that $\mathcal{S}$ – which by the above is a locally finitely generated submodule – is stable under time derivatives in this sense: for all smooth 1-parameter families $\{\boldsymbol{\alpha}_\lambda\} \subset \Gamma_c(A)$ such that $\boldsymbol{\alpha}_\lambda \in \mathcal{S}$ for all $\lambda$ and $\boldsymbol{\alpha}_0 = 0$, the time derivative $\frac{d}{d\lambda}|_{\lambda=0}\boldsymbol{\alpha}_\lambda$ lies in $\mathcal{S}$. Then $\mathcal{S}$ is involutive.

# 4 Global Integration of a singular subalgebroid

In this section we introduce a global notion of integration for a singular subalgebroid $\mathcal{B}$ (Def. 4.1), which we call *integral*, and prove that the holonomy groupoid $H^G(\mathcal{B})$ is an integral (Thm. 4.19). We also show that $\mathcal{B}$ admits a smooth integral iff it is projective (Cor. 4.14). In that case $H^G(\mathcal{B})$ is the minimal smooth integral, and all other smooth integrals are coverings of $H^G(\mathcal{B})$ (Prop. 4.16 and 4.17). When $\mathcal{B}$ is a Lie subalgebra of a Lie algebra, all integrals are smooth (Prop. 4.7).

## 4.1 Definition of integral

We now define the notion of integral of a singular subalgebroid $\mathcal{B}$, as a diffeological groupoid that differentiates to $\mathcal{B}$ (in the sense of §3) and satisfies two additional properties. Definition 4.1 below is inspired by [1, Appendix A].

**Definition 4.1 (Integrals).** Let $G$ be a Lie groupoid over $M$, $H$ be a diffeological groupoid over $M$, and $\Psi\colon H \to G$ a smooth morphism of diffeological groupoids covering $Id_M$.

Let $\mathcal{B}$ be a singular subalgebroid of $A := AG$. We say that $(H, \Psi)$ is an **integral of $\mathcal{B}$ over $G$** if:

a) $(H, \Psi)$ differentiates to $\mathcal{B}$ (Def. 3.3).
(In particular $H$ is holonomy-like, see Def. 2.12, and Assumption 3.1 is satisfied.)

b) The diffeology of $H$ is generated by open maps (Def. 2.18).

c) The morphism $\Psi$ is *coarsely injective*, in the following sense: There exists a neighborhood $\overset{\circ}{H}$ of the identity $1_M$ in $H$ such that if $\mathbf{b} \subset \overset{\circ}{H}$ is a bisection (see Def. 2.7) carrying the identity bisection of $G$, then $\mathbf{b} \subset 1_M$.

Several remarks are in order.

**Remark 4.2** (The three conditions)**.** We comment on the three conditions appearing in Def. 4.1.

a) Condition a) is certainly expected for any reasonable notion of integration.

b) Condition b) implies that $H$ is source-submersive – i.e. through every point there passes a bisection –, thanks to Prop. 2.17 i). This makes condition c) a meaningful one.

   The fact that condition b) is sufficient to guarantee that $H$ is source-submersive, and the fact that checking when a map is open is generally easier than checking when it is a local subduction, are the main reasons why in Definition 4.1 we chose to use condition b) rather than the stronger condition that the diffeology of $H$ be generated by local subductions (Def. 2.21).

c) The coarsely injective requirement on $\Psi$ is imposed to ensure that $H$ is "not larger than necessary". For example, when $M$ is a point, so that $H$ and $G$ are groups, $\Psi$ must be an injection in a neighborhood of the identity element of $H$.

**Remark 4.3** (Conditions on $H$)**.** In the definition of integral $(H, \Psi)$, there are two requirements on the diffeological groupoid $H$: that it be holonomy-like (Def. 2.12), and the openness condition b). Each of these two requirements consist of the existence of a generating set of plots with certain properties, in some neighborhood $\overset{\circ}{H}$ of $1_M$ in $H$ (see Lemma 2.20).

i) If the generating set of plots $\{\chi_i : \mathcal{O}_{\chi_i} \to \overset{\circ}{H}\}$ required by the holonomy-like property consists of open maps, then both requirements are satisfied. A special case is when the diffeology of $H$ is generated by local subductions (Def. 2.21), according to Remark 2.23. Spelled out, this means that for each $i$ the image $\chi_i(\mathcal{O}_{\chi_i})$ is an open subset of $H$ and $\chi_i : \mathcal{O}_{\chi_i} \to \chi_i(\mathcal{O}_{\chi_i})$ is a local subduction, by Lemma 2.22. This is exactly what occurs for the examples of integrals we display in §4.3 and §4.4, namely smooth integrals and the holonomy groupoid.

ii) The above two requirements together, by Lemma 1.18, imply that the generating set $\{\chi_i : \mathcal{O}_{\chi_i} \to \overset{\circ}{H}\}$ for the holonomy-like property can be chosen so that the image of each map $\chi_i$ is an open subset (but not necessarily that the $\chi_i$ are open maps). This observation supports the scenario outlined in item i).

**Remark 4.4** (Restatement of Assumption 3.1)**.** In Def. 4.1, Assumption 3.1 admits an equivalent characterization as Assumption B.3 on the liftability of 1-parameter groups of global bisections of $G$. We prove this in Corollary B.2.

This equivalent characterization will be used only in the proof of Prop. 5.8, which states that the image $\Psi(H)$ is the same for all integrals $(H, \Psi)$ of $\mathcal{B}$.

**Remark 4.5** (Functoriality)**.** The differentiation functor introduced in §A.4 is well-behaved when restricted to integrals. Indeed Lemma A.10 applies automatically to integrals, thanks to Remark 4.2 b)

**Remark 4.6.** Notice that the conditions for $(H, \Psi)$ to be an integral are all conditions on a neighbourhood of the units $1_M$ of $H$: indeed, condition a) on differentiation (in particular $H$ being holonomy-like) and condition c) on being coarsely injective are in terms of such a neighbourhood, and condition b) also is thanks to Lemma 2.20. Recall that this is a well known feature of classical Lie groupoids (in the smooth category): The Lie algebroid of a Lie groupoid is determined by a neighborhood of the identity.

## 4.2 The integrals of Lie subalgebras

Here we consider the special case of a Lie subalgebra of a (finite dimensional, real) Lie algebra. We show that all integrals are actually Lie groups.

**Proposition 4.7.** Let $\mathfrak{g}$ a Lie algebra, $\mathfrak{k}$ a Lie subalgebra, fix a Lie group $G$ integrating $\mathfrak{g}$. Let $(H, \Psi)$ be an integral of $\mathfrak{k}$. Then there exists a Lie group structure whose underlying diffeological group is $H$.

**Remark 4.8.** The above Lie group structure, which we denote by $H_{\text{Lie}}$, is unique by Ex. 1.15 a). Further $H_{\text{Lie}}$ integrates the Lie algebra $\mathfrak{k}$, and $\Psi$ – viewed as a map on $H_{\text{Lie}}$ – integrates the inclusion $\mathfrak{k} \hookrightarrow \mathfrak{g}$, see the proof of Prop. 4.7 below. Thus $H_{\text{Lie}}$ is a covering of the (unique) connected Lie subgroup of $G$ integrating the Lie subalgebra $\mathfrak{k}$.

We recall what it means that $(H, \Psi)$ is an integral of $\mathfrak{k}$ (Def. 4.1):

a1) $H$ is a diffeological group which is holonomy-like (we will not use the latter condition).

b) The diffeology of $H$ is generated by open maps (we will not use this condition).

a2) $\Psi \colon H \to G$ is a smooth morphism of diffeological groups such that

$$\begin{aligned}\mathfrak{k} &= \left\{ \frac{d}{d\lambda}|_{\lambda=0}(\Psi \circ \mathbf{b}_\lambda) \colon \ \{\mathbf{b}_\lambda\}_{\lambda \in I} \text{ family of elements of } H \text{ s.t. } \mathbf{b}_0 = Id_M \right\} \\ &= \left\{ \frac{d}{d\lambda}|_{\lambda=0}(\Psi \circ \mathbf{b}_\lambda) \colon \ \{\mathbf{b}_\lambda\}_{\lambda \in I} \text{ 1-parameter group of elements of } H \right\}\end{aligned}$$

c) $\Psi$ is injective on a neighborhood of the unit in $H$.

*Proof of Prop. 4.7.* First observe that if $\{\mathbf{b}_\lambda\}$ is a 1-parameter group in $H$, then $\{\Psi \circ \mathbf{b}_\lambda\}$ is a 1-parameter group in the Lie group $G$, hence

$$\Psi \circ \mathbf{b}_\lambda = \exp(\lambda v) \tag{4.1}$$

where $v := \frac{d}{d\lambda}|_{\lambda=0}(\Psi \circ \mathbf{b}_\lambda)$ and exp is the Lie group exponential map of $G$.
Denote by $K$ the unique connected Lie subgroup of $G$ integrating the Lie subalgebra $\mathfrak{k}$.
Claim: $\Psi(H) = K$.
*Proof of claim:*

"$\supset$" There is a neighborhood of the unit in $K$ consisting of elements of the form $\exp(v)$ where $v \in \mathfrak{k}$. By assumption, there is a 1-parameter group $\{\mathbf{b}_\lambda\}$ in $H$ such that $v = \frac{d}{d\lambda}|_{\lambda=0}(\Psi \circ \mathbf{b}_\lambda)$. By the above observation, taking $\lambda = 1$, we get $\Psi \circ \mathbf{b}_1 = \exp(v)$. Conclude using that $K$, being connected, is generated by any neighborhood of the unit.

"$\subset$" For any diffeological space, the set of points that can be connected to a given point by a smooth path is open in the $D$-topology. (This can be proven exactly as [18, Lemma 1.8], which states that the $D$-topology is locally arc-connected. Notice that while the composition of two smooth paths is not smooth in general, it is if we first apply a suitable time reparametrization.) In particular[8], the set of points in $H$ that can be connected to the unit $\mathbf{e}$ by a smooth path is an open subset $\breve{H}$ of $H$.

Let $h \in \breve{H}$ and $\{\mathbf{b}_\lambda\}_{\lambda\in[0,1]}$ be a smooth path in $H$ from the identity element $\mathbf{e}$ to $h$. For every fixed $\lambda_0 \in (0,1)$, consider the smooth path $\{\mathbf{b}_{\mu+\lambda_0} \cdot \mathbf{b}_{\lambda_0}^{-1}\}_{\mu\in I}$ going through $\mathbf{e}$ at time $\mu = 0$, where $I$ is an open interval containing zero. Since $H$ differentiates to $\mathfrak{k}$, we have

$$\frac{d}{d\mu}|_{\mu=0}\left(\Psi(\mathbf{b}_{\mu+\lambda_0} \cdot \mathbf{b}_{\lambda_0}^{-1})\right) \in \mathfrak{k}.$$

Since $\Psi$ is a group morphism, the left hand side equals

$$\frac{d}{d\mu}|_{\mu=0}R_{\Psi(\mathbf{b}_{\lambda_0}^{-1})}(\Psi(\mathbf{b}_{\mu+\lambda_0})) = (R_{\Psi(\mathbf{b}_{\lambda_0}^{-1})})_*\frac{d}{d\lambda}|_{\lambda=\lambda_0}(\Psi(\mathbf{b}_\lambda))$$

where $R$ denotes right-translation on $G$. Hence $\frac{d}{d\lambda}|_{\lambda=\lambda_0}(\Psi(\mathbf{b}_\lambda))$ lies in the right-invariant distribution $\vec{\mathfrak{k}}$. Repeating for all $\lambda_0 \in (0,1)$ we see that the curve $[0,1] \ni \lambda \mapsto \Psi(\mathbf{b}_\lambda)$ lies in the leaf of $\vec{\mathfrak{k}}$ through $\mathbf{e}$, which is exactly the Lie subgroup $K$. We conclude that $\Psi(h) \in K$. △

The morphism of topological groups $\Psi\colon H \to G$ is injective in a neighborhood of the unit, and by the claim its image is the Lie group $K$. Hence $K$, via $\Psi$, induces a Lie group structure on $H$, which we denote by $H_{\mathrm{Lie}}$, and whose Lie algebra is isomorphic to $\mathfrak{k}$. We are left with showing that the diffeological structure underlying the manifold structure of $H_{\mathrm{Lie}}$ coincides with the diffeological structure of $H$.

Claim: The identity map is an isomorphism of diffeological groups between $H$ and $H_{\mathrm{Lie}}$.

*Proof of claim:* We proceed in two steps.

The identity $H \to H_{\mathrm{Lie}}$ is smooth: Since $\Psi\colon H \to G$ is smooth, the induced surjective map to the Lie group $K$ is a smooth map: for any plot $\chi$ into $H$, the composition $\Psi \circ \chi$ is a smooth plot into $K$, since $K$ is an initial submanifold of $G$ [20, Thm. 19.25].

The identity $H_{\mathrm{Lie}} \to H$ is a smooth: Thanks to Lemma 4.10 below, it suffices to show that all plots for the *strong* diffeology of $H_{\mathrm{Lie}}$ are also plots for the diffeology on $H$. That is, we have to show that all rays in $H_{\mathrm{Lie}}$ – which, as is well-known, are of the form $\lambda \mapsto \exp_{H_{\mathrm{Lie}}}(\lambda v)$ for elements $v$ in the Lie algebra $T_{\mathbf{e}}H_{\mathrm{Lie}}$ – are plots for $H$. To show this, fix $v \in T_{\mathbf{e}}H_{\mathrm{Lie}}$, and consider $\Psi_* v \in \mathfrak{k}$. By assumption, there is a (smooth) 1-parameter group $\{\mathbf{b}_\lambda\}_{\lambda\in I}$ in $H$ such that $\Psi_* v = \frac{d}{d\lambda}|_{\lambda=0}(\Psi \circ \mathbf{b}_\lambda)$, and by the above observation $\Psi \circ \mathbf{b}_\lambda = \exp_K(\lambda \Psi_* v)$. By the injectivity of $\Psi$ near the unit, this implies that $\mathbf{b}_\lambda = \exp_{H_{\mathrm{Lie}}}(\lambda v)$ for values of $\lambda$ close enough to zero, showing that $\lambda \mapsto \exp_{H_{\mathrm{Lie}}}(\lambda v)$ is smooth for the diffeology of $H$. □

**Remark 4.9.** In the setting of Prop. 4.7, recall that $K$ is the holonomy groupoid $H^G(\mathfrak{k})$ of $\mathfrak{k}$, by item d) in Ex. 1.10. The integral $H$ quotients to $H^G(\mathfrak{k})$ (see the proof of Prop. 4.7); this fact will be generalized in Prop. 4.17.

We recall a few notions from [27, §6]. Given a diffeological group $H$, a **ray** is a smooth group homomorphism $\mathbb{R} \to H$. The **strong diffeology** on $H$ is the smallest diffeology on $H$ containing all rays and so that the group multiplication and inversion are smooth. All plots for the strong diffeology are plots for the original diffeology on $H$.

[8] Alternatively, this also follows from Condition b) in Def. 4.1.

For Lie groups the converse holds. This is the content of the following result, which was used in the proof of Prop. 4.7 above and was stated without proof by Souriau [27, §6.13]. We sketch a proof for the sake of completeness.

**Lemma 4.10.** Let $H$ be a Lie group. Then the strong diffeology of $H$ coincides with the Lie group diffeology.

*Proof.* It is well-known that the Lie group exponential map $\exp\colon \mathfrak{h} \to H$, once restricted to a suitable open neighborhood of the origin, provides a chart for the manifold structure on the Lie group $H$. In view of the text just before this lemma, we have to show that all plots for the manifold diffeology on $H$ are plots for the strong diffeology. It suffices to show that exp is smooth with respect to the strong diffeology, by Ex. 1.15 a).

To this aim, let $n := \dim(H)$, and fix a basis $v_1, \dots, v_n$ of the Lie algebra $\mathfrak{h}$. For a small enough neighborhood $U \subset \mathbb{R}^n$ of the origin, the map

$$S\colon U \to H,\ \ (\lambda_1, \dots, \lambda_n) \mapsto \exp(\lambda_1 v_1) \cdot \exp(\lambda_2 v_2) \cdot \dots \cdot \exp(\lambda_n v_n)$$

is a smooth map to the Lie group $H$ (actually, a diffeomorphism onto its image), and it is also smooth for the strong diffeology, since it is obtained multiplying rays. Consider the map

$$B\colon U \to \mathfrak{h},\ \ (\lambda_1, \dots, \lambda_n) \mapsto BCH\Big(\lambda_1 v_1, BCH\big(\dots, BCH(\lambda_{n-1} v_{n-1}, \lambda_n v_n)\big)\Big),$$

where “BCH” denotes the Baker-Campbell-Hausdorff formula. The relation between these two maps is $S = \exp \circ B$. Thus $B = \exp^{-1} \circ S\colon U \to \mathfrak{h}$ is a diffeomorphism onto its image. Restricting to this image we have

$$\exp = S \circ B^{-1}.$$

This shows that exp is a plot for the strong diffeology as the composition the smooth map $B^{-1}$ between open subsets of vector spaces and the plot $S$ for the strong diffeology. □

For future reference, we present a proof of Prop. 4.7 which does not use Assumption 3.1, but instead uses that $H$ is holonomy-like[9] .

*Alternative proof of Prop. 4.7.* We modify the proof above in two places, for which we provide an alternative argument.

$\underline{\Psi(H) \supset K}$: Fix a generating set for the diffeology of a neighborhood $\overset{\circ}{H}$ of the unit element $\mathbf{e}$ as in Def. 2.12. In this generating set, fix a plot $(\chi, \mathcal{O}_\chi)$ containing a point $e$ such that $\chi(e) = \mathbf{e}$.

Claim: The derivative of $\Psi \circ \chi\colon \mathcal{O}_\chi \to K$ at the point $e$ is surjective.

*Proof of claim:* Notice that the above map is well-defined since $\Psi(H) \subset K$. Take a vector $w \in \mathfrak{k}$. There is a smooth path $\mathbf{b}_\lambda$ in $H$ with $\mathbf{b}_0 = \mathbf{e}$ and $\frac{d}{d\lambda}|_{\lambda=0}(\Psi \circ \mathbf{b}_\lambda) = w$, since $H$ differentiates to $\mathfrak{k}$.

[9]By Def. 2.12 this means that there exists an open neighborhood $\overset{\circ}{H}$ of the unit 1 in $H$ and a generating set of plots $\mathcal{X} = \{\chi : \mathcal{O}_\chi \to \overset{\circ}{H}\}$ in $\mathcal{P}(\overset{\circ}{H})$ such that: given any plots $\chi, \chi'$ in $\mathcal{X}$ and points $e \in \mathcal{O}_\chi, e' \in \mathcal{O}_{\chi'}$ such that $\chi(e) = \chi'(e') = 1$, there exists a smooth map $k\colon \mathcal{O}_{\chi'} \to \mathcal{O}_\chi$ with $k(e') = e$ and $\chi \circ k = \chi'$, after shrinking $\mathcal{O}_{\chi'}$ to a smaller neighborhood of $e'$ if necessary.

There is a plot $(\chi', \mathcal{O}_{\chi'})$ in the generating set to which $\mathbf{b}_\lambda$ lifts, in the sense that there is a smooth path $\gamma'_\lambda$ in $\mathcal{O}_{\chi'}$, possibly defined on a smaller interval, such that $\chi' \circ \gamma'_\lambda = \mathbf{b}_\lambda$. Denote $e' := \gamma'_0$.

$$\begin{array}{ccc} H & \xrightarrow{\Psi} & K \\ {\scriptstyle \chi'}\big\uparrow & & \\ \mathcal{O}_{\chi'} & & \end{array}$$

By the holonomy-like property, there exists a smooth map $k\colon \mathcal{O}_{\chi'} \to \mathcal{O}_\chi$ with $k(e') = e$ and $\chi \circ k = \chi'$, after shrinking $\mathcal{O}_{\chi'}$ to a smaller neighborhood of $e'$ if necessary. So $\gamma_\lambda := k \circ \gamma'_\lambda$ is a curve in $\mathcal{O}_\chi$ which goes through $e$ at time $\lambda = 0$, and such that $\chi \circ \gamma_\lambda = \mathbf{b}_\lambda$. In particular, $\frac{d}{d\lambda}|_{\lambda=0}((\Psi \circ \chi) \circ \gamma_\lambda) = w$, proving the claim. △

Nearby $e$, by the claim, $\Psi \circ \chi$ is a submersion, thus an open map, and hence its image $(\Psi \circ \chi)(\mathcal{O}_\chi)$ contains an open neighborhood of $\mathbf{e}$ in the Lie group $K$. Therefore the group $\Psi(H)$ contains an open neighborhood of $\mathbf{e}$ in $K$, and consequently it contains $K$.

<u>The identity $H_{\text{Lie}} \to H$ is a smooth.</u>

As in the proof of Prop. 4.7, it suffices to show that all rays in $H_{\text{Lie}}$ are plots for $H$. Let $v$ be a vector in the Lie algebra of $H_{\text{Lie}}$. We need to show that $\mathbb{R} \ni \lambda \mapsto \exp_{H_{\text{Lie}}}(\lambda v)$ is a plot for $H$. By the previous claim, $\Psi \circ \chi\colon \mathcal{O}_\chi \to K$ is a submersion nearby $\mathbf{e}$, so we can lift the path $\exp_K(\lambda \Psi_* v)|_I$ in $K$ to a smooth path $\tau_\lambda$ in $\mathcal{O}_\chi$ with $\tau_0 = \mathbf{e}$, for an open interval $I$ around zero. Hence $\chi \circ \tau_\lambda$ is a smooth path in $H$ with two properties: it goes through $1_H$ at time $\lambda = 0$, an it maps to $\exp_K(\lambda \Psi_* v)|_I$ under $\Psi$. On the other hand, $\exp_{H_{\text{Lie}}}(\lambda v)|_I$ is also a path in $H$ with the same two properties. Since $\Psi$ is injective in a neighborhood of $\mathbf{e}$ (w.r.t. the D-topology on $H$), it follows that the two paths must agree (after possibly shrinking $I$), i.e.

$$(\chi \circ \tau_\lambda)|_I = \exp_{H_{\text{Lie}}}(\lambda v)|_I.$$

Now, the curve on the left is a *smooth* curve into $H$, thus the one on the right too. Using the fact that left-multiplication by $\exp_{H_{\text{Lie}}}(\lambda_0 v)$ is a smooth map $H \to H$ for all $\lambda_0 \in \mathbb{R}$, we conclude that $\mathbb{R} \to H, \lambda \mapsto \exp_{H_{\text{Lie}}}(\lambda v)$ is a smooth curve. □

## 4.3 Smooth integrals

In this subsection we determine which singular subalgebroids admit integrals (in the sense of definition 4.1) which are smooth. For any such singular subalgebroid, we exhibit all the smooth integrals.

Let $H \rightrightarrows M$ be a diffeological groupoid, in the sense of Def. 2.1. If there exists a manifold structure on $H$ making the above diffeological groupoid a Lie groupoid, this manifold structure is unique, by Ex. 1.15 a).

**Definition 4.11.** An integral $(H, \Psi)$ of a singular subalgebroid $\mathcal{B}$ is said to be a **smooth integral** if the diffeological groupoid $H$ admits a (necessarily unique) manifold structure making it a Lie groupoid.

Notice that if $(H, \Psi)$ is a smooth integral of $\mathcal{B}$, then $\Psi$ is automatically a Lie groupoid morphism, by Ex. 1.15 a). Further, $\mathcal{B} = \Psi_*(\Gamma_c(AH))$ by Def. 4.1 a).

**Lemma 4.12.** Assume $H$ is a Lie groupoid and $\Psi\colon H \to G$ a Lie groupoid morphism covering $Id_M$. Define $\mathcal{B} := \Psi_*(\Gamma_c(AH))$. Then $\Psi\colon H \to G$ is an integral for $\mathcal{B}$ iff the Lie algebroid morphism $\Psi_*\colon AH \to AG$ is almost injective[10].

*Proof.* "⇒" Let $\boldsymbol{\alpha} \in \Gamma_c(AH)$ such that $\Psi_*\boldsymbol{\alpha} = 0$. It suffices to show that $\boldsymbol{\alpha} = 0$. Define a smooth family of bisections $\mathbf{b}_\lambda\colon M \to H$ by $\mathbf{b}_\lambda(x) = exp_x(\lambda\overrightarrow{\boldsymbol{\alpha}})$. For all $x \in M$,

$$\Psi \circ exp_x(\lambda\overrightarrow{\boldsymbol{\alpha}}) = exp_x(\lambda\Psi_*\overrightarrow{\boldsymbol{\alpha}}) = x,$$

where in the first equality we used that $\Psi$ is a Lie groupoid morphism and in the second that $\Psi_*\overrightarrow{\boldsymbol{\alpha}} = \overrightarrow{\Psi_*\boldsymbol{\alpha}}$ is the zero vector field on $G$. This shows that $\Psi \circ \mathbf{b}_\lambda$ is the identity bisection of $G$, for all $\lambda$ sufficiently close to zero, hence by a) in Def. 4.1 we obtain that $\mathbf{b}_\lambda$ is the identity section of $H$ for all such $\lambda$, and therefore that $\boldsymbol{\alpha} = 0$.

"⇐" Condition a) in Def. 4.1 is checked in Prop. A.6. Condition b) holds for all Lie groupoids, see Remark 2.26 and Remark 2.23. Now we check condition c). Since $\Psi_*$ is almost injective, for all $x$ lying in a dense subset of $M$ we have that the restriction of $\Psi$ to the fiber $\mathbf{s}_H^{-1}(x)$ is injective nearby $x$. This implies that $\Psi$ is injective for bisections nearby the identity. □

By Lemma 4.12, if $\Psi\colon H \to G$ a smooth integral of a singular subalgebroid $\mathcal{B}$, then $\mathcal{B}$ is necessarily a projective singular subalgebroid [6], whose underlying Lie algebroid is isomorphic to $AH$.

Conversely, every projective singular subalgebroid admits a smooth integral, as the following example shows.

**Example 4.13.** Let $\mathcal{B}$ be a projective singular subalgebroid of $A = AG$, hence there exists a Lie algebroid $B$ with an almost injective morphism

$$\psi\colon B \to A$$

inducing an isomorphism $\Gamma_c(B) \cong \mathcal{B}$. Let $H$ be a Lie groupoid integrating $B$ with a morphism of Lie groupoids $\Psi\colon H \to G$ integrating $\psi$. Then $(H, \Psi)$ is a smooth integral of $\mathcal{B}$ over $G$, as follows immediately from Lemma 4.12. We point out that $H$ as above always exists by [6]; as $H$ we can take the source simply connected Lie groupoid integrating $B$, and by [6] we can also take the holonomy groupoid $H^G(\mathcal{B})$. Both of them are therefore smooth integrals.

As an immediate corollary of the discussion above, we determine the singular subalgebroids for which there exists a smooth integral.

**Corollary 4.14.** Let $\mathcal{B}$ be a singular subalgebroid of $AG$. Then $\mathcal{B}$ admits a smooth integral iff $\mathcal{B}$ is projective. Further, the holonomy groupoid $H^G(\mathcal{B})$ is a smooth integral.

**Remark 4.15.** Given a projective singular subalgebroid, in general we do not know whether all its integrals are smooth. We do not even know this for the subclass of singular subalgebroid arising from wide Lie subalgebroids, except for the simple case of Lie subalgebras (see Prop. 4.7).

When $\mathcal{B}$ is projective, the holonomy groupoid $H^G(\mathcal{B})$ is minimal among all of its possible *smooth* integrals. Note that non-smooth integrals of a projective singular subalgebroid $\mathcal{B}$ may exist, for instance the graph of a regular foliation. (An example is the open Möbius band with a foliation by circles, see e.g. [25, p. 154]). Recall that the morphism $\Phi$ was defined in §1.3.

[10] Recall that this means that $\Psi_*$ fiber-wise injective on a dense subset of $M$.

**Proposition 4.16** (Minimality among smooth integrals)**.** Let $\mathcal{B}$ be a singular subalgebroid of $AG$. Let $(H, \Psi)$ be a *smooth* integral of $\mathcal{B}$ over $G$. Then there is a surjective morphism of Lie groupoids $H \to H^G(\mathcal{B})$ such that this diagram commutes:

$$\begin{array}{ccc} H & & \\ \downarrow & \searrow^{\Psi} & \\ H^G(\mathcal{B}) & \xrightarrow{\Phi} & G \end{array}$$

*Proof.* As a consequence of Lemma 4.12, $\mathcal{B} = \Psi_*(\Gamma_c(AH))$ is projective. Since $\Psi : H \to G$ is a Lie groupoid morphism, it is a bisubmersion for $\mathcal{B}$ by item a) in Ex. 1.10. Then item b) in Ex. 1.10 implies the statement. (The smoothness of the quotient map $H \to H^G(\mathcal{B})$ follows from $\Psi$ being a bisubmersion for $\mathcal{B}$). □

The following proposition characterizes the smooth integrals. It is the analogue, for projective singular subalgebroids, of a theorem of Moerdijk-Mrčun on the integration of wide Lie subalgebroids [24, Thm. 2.3].

**Proposition 4.17** (Characterization of smooth integrals)**.** Let $\mathcal{B}$ be a projective singular subalgebroid of $AG$. Denote by $B$ the corresponding Lie algebroid, satisfying $\Gamma_c(B) \cong \mathcal{B}$. The smooth integrals of $\mathcal{B}$ over $G$ are precisely the pairs

$$(H, \Phi \circ \pi),$$

where $\pi \colon H \to H^G(\mathcal{B})$ is a surjective Lie groupoid morphism integrating $Id_B$ (in particular $H$ is a Lie groupoid integrating $B$).

*Proof.* The pairs $(H, \Phi \circ \pi)$ in the statement are smooth integrals of $\mathcal{B}$ over $G$, as explained in Ex. 4.13. All smooth integrals are of this kind, by Prop. 4.16. □

**Remark 4.18.** When $\mathcal{B} = \Gamma_c(B)$ for a wide Lie subalgebroid $B \subset AG$, its smooth integrals recover exactly the Lie groupoid morphisms of [24, Thm. 2.3].

### 4.4 The holonomy groupoid is an integral

Here we prove that $H^G(\mathcal{B})$ satisfies definition 4.1.

**Theorem 4.19.** The holonomy groupoid $(H^G(\mathcal{B}), \Phi)$ is an integral of $\mathcal{B}$ over $G$

*Proof.* We already proved in Thm. A.7 that $(H^G(\mathcal{B}), \Phi)$ differentiates to $\mathcal{B}$. Further, $H^G(\mathcal{B})$ satisfies Condition b) of definition 4.1, by Prop. 1.9 and Examples 1.15 d). We now prove that $(H^G(\mathcal{B}), \Phi)$ satisfies the injectivity condition c) of definition 4.1. Let $\mathbf{b}$ be a bisection of $H^G(\mathcal{B})$ carrying the identity bisection of $G$. Then locally there exists a bisubmersion $(U, \varphi, G)$ with a bisection $\mathbf{u}$ such that $\mathbf{b}$ is the image of $\mathbf{u}$ by the quotient map $q_U : U \to H^G(\mathcal{B})$. So $\mathbf{u}$ also carries the identity bisection of $G$, which means that $q_U(\mathbf{u}) \subset 1_M$ by the second bullet point in §1.3 . Whence $\mathbf{b} \subset 1_M$.

□

**Remark 4.20.** In view of Theorem 4.19, the question arises of whether the groupoid $H^G(\mathcal{B})$ is a *minimal* integral of $\mathcal{B}$ over $G$, in some sense. In proposition 4.16 we gave a partial answer for projective singular subalgebroids, in the spirit of [24]. In Cor. 5.7 we show that for a specific smooth holonomy groupoid, any of its non-trivial quotients fails to be an integral.

## 5 The graph of a singular subalgebroid

It is well known that, given a regular foliation, one can attach to it the holonomy groupoid as well as the graph of the foliation, namely the groupoid defined by the equivalence relation of "belonging to the same leaf". The graph generally fails to be a smooth manifold, and this is one of the *raison d'être* of the holonomy groupoid: the holonomy groupoid is the minimal Lie groupoid[11] with Lie algebroid isomorphic to the tangent distribution of the foliation.

In §5.1, for any singular subalgebroid $\mathcal{B}$, we show that there is a natural diffeology on the graph, called path holonomy diffeology, making it into a diffeological groupoid. We remark that the construction of the path holonomy diffeology makes essential use of the holonomy groupoid. We show that the graph differentiates to $\mathcal{B}$ (Prop. 5.5). In particular, the path holonomy diffeology captures enough information to recover $\mathcal{B}$. However the graph is *not* an integral of $\mathcal{B}$ in general, since it generally violates Condition b) in the definition of integral (Prop. 5.6). This raises the question of whether the holonomy groupoid $H^G(\mathcal{B})$ is "minimal" – in a sense to make precise – among all integrals of $\mathcal{B}$ over $G$, see Remark 4.20. In §5.2 we relate the graph and the path holonomy diffeology to arbitrary integrals of $\mathcal{B}$. In particular, we show that the image of any integral of $\mathcal{B}$ under its canonical map to $G$ is exactly the graph (Prop. 5.8).

Let us start giving the definitions.

**Definition 5.1.** Let $\mathcal{B}$ be a singular subalgebroid of $AG$. The **graph of $\mathcal{B}$ over** $G$ is the subgroupoid $R^G(\mathcal{B})$ of $G$ defined as the image of the groupoid morphism $\Phi : H^G(\mathcal{B}) \to G$.

In the case of a singular foliation $(M, \mathcal{F})$, the morphism $\Phi$ is the target-source map, so the graph $R(\mathcal{F}) := R^{M\times M}(\mathcal{F})$ is the subgroupoid of the pair groupoid $M \times M$ given by the equivalence relation of "belonging to the same leaf of $\mathcal{F}$".

We now give an explicit description of the graph of an arbitrary singular subalgebroid $\mathcal{B}$ of a Lie algebroid $AG$. This description makes clear that the graph does not depend on $\mathcal{B}$ but just on its image under the evaluation map.

**Remark 5.2** (Description of the graph)**.** For every leaf $L$ of (the singular foliation induced by) $\mathcal{B}$, there is transitive Lie subalgebroid $B_L \to L$ of $AG$, obtained evaluating point-wise the elements of $\mathcal{B}$, see [6]. By right-translation we obtain a regular foliation $\overrightarrow{B_L}$ on $G_L = \mathbf{t}^{-1}(L)$, see also [13, Prop 5.1] for a different description. (When the leaf $L$ is immersed but not embedded, $G_L$ is endowed with the differentiable structure coming from the canonical bijection $G_L \cong L \times_{(\iota,\mathbf{t})} G$, where $\iota : L \to M$ is the inclusion.) The union of the leaves of $\overrightarrow{B_L}$ which intersect the manifold of unit elements is a subgroupoid of $G_L$, which we denote by $R^G(\mathcal{B})_L$. The graph agrees with the union

$$R^G(\mathcal{B}) := \bigcup_L R^G(\mathcal{B})_L,$$

where $L$ ranges through all leaves of $\mathcal{B}$.

This follows from the following facts, shown in [6]: i) denoting by $\mathcal{B}_L$ the transitive Lie algebroid over $L$ whose sections are given by by $\mathcal{B}/I_L\mathcal{B}$, we have an evaluation map $ev_L \colon \mathcal{B}_L \to A$, whose image is $B_L$; ii) the Lie algebroid morphism $ev_L$ is integrated by the restriction to $H^G(\mathcal{B})|_L$ of $\Phi : H^G(\mathcal{B}) \to G$.

[11] Although the holonomy groupoid is not always a Hausdorff manifold.

## 5.1 The graph as diffeological groupoid

The holonomy groupoid, via the canonical morphism $\Phi\colon H^G(\mathcal{B}) \to G$, induces a diffeology on the graph, which we call path-holonomy diffeology. The graph $R^G(\mathcal{B})$ is thus naturally a diffeological groupoid, and it comes with the inclusion map to $G$. In Prop. 5.5 we show that it differentiates to $\mathcal{B}$. However it is not an integral of $\mathcal{B}$ in general, since it generally violates Condition b) in the definition of integral (Prop. 5.6).

### 5.1.1 The subspace diffeology on the graph

The **subspace diffeology** on the the graph is the one obtained restricting the diffeology of $G$ arising from its smooth structure. The plots in this diffeology are all the smooth maps $f : U \subseteq \mathbb{R}^n \to G$ such that $f(U) \subseteq \Phi(H^G(\mathcal{B}))$. Due to the following reasons, the subspace diffeology captures little information and is thus not useful for us:

- This diffeology does not depend on the map $\Phi$ but rather on its image.
- This diffeology is not holonomy-like in general (see Example 5.3 b) below).
- $R^G(\mathcal{B})$, with the subspace diffeology, may *not* differentiate[12] to $\mathcal{B}$. This follows from Examples 5.3 below.

**Example 5.3.** a) Take $M = \mathbb{R}$ with the foliation $\mathcal{F}_i = span_{C_c^\infty(\mathbb{R})}\langle x^i \partial_x \rangle$, where $i = 1, 2$. Then $R(\mathcal{F}_1) = R(\mathcal{F}_2) = (\mathbb{R}^+ \times \mathbb{R}^+) \cup (\mathbb{R}^- \times \mathbb{R}^-) \cup \{(0, 0\}$. But the differentiation of $R(\mathcal{F}_i)$ with respect to the $(\mathbb{R} \times \mathbb{R})$-diffeology gives the module generated by all vector fields on $\mathbb{R}$ which vanish at 0.

b) Take $M = \mathbb{R}$ with the foliation $\mathcal{F} = span_{C_c^\infty(\mathbb{R})}\langle f \partial_x \rangle$, where $f$ is a smooth function vanishing on $\mathbb{R}^-$ and which is no-where vanishing on $\mathbb{R}^+$. Then $R(\mathcal{F}) = \bigcup_{x \leqslant 0} \{(x, x)\} \cup (\mathbb{R}^+ \times \mathbb{R}^+)$. Its differentiation with respect to the $(\mathbb{R} \times \mathbb{R})$-diffeology never returns $\mathcal{F}$, indeed it gives the module generated by all vector fields on $\mathbb{R}$ which vanish on $\mathbb{R}^-$. But this module is not locally finitely generated, so this diffeology cannot be holonomy-like.

It is not surprising that $R^G(\mathcal{B})$ with the subspace diffeology does not differentiate to $\mathcal{B}$ in general. In the case of singular foliations, there may be more than one module of vector fields which induce the same partition into leaves. However, for a singular foliation $(M, \mathcal{F})$, the diffeological groupoid $R^{M\times M}(\mathcal{F})$ with the $M \times M$-diffeology depends only on the partition into leaves it induces on $M$, hence it does not contain enough information to recover $\mathcal{F}$ by differentiation.

### 5.1.2 The path holonomy diffeology on the graph

Much more useful is the **path-holonomy diffeology**, i.e. the quotient diffeology induced by the surjective map $\Phi\colon H^G(\mathcal{B}) \to \Phi(H^G(\mathcal{B}))$. This diffeology is generated by the plots $\Phi \circ q_U$, where $(U, \varphi, G)$ is a bisubmersion in a path-holonomy atlas of $\mathcal{B}$. With this diffeology, the inclusion $\iota\colon R^G(\mathcal{B}) \hookrightarrow G$ is a smooth morphism, since $\Phi \circ q_U\colon U \to G$ is smooth. The path-holonomy diffeology makes $R^G(\mathcal{B})$ into a diffeological groupoid.

The following two lemmas show that the path-holonomy diffeology is well-behaved.

[12] Here we are applying Def. 3.3 despite the fact that $R^G(\mathcal{B})$ with the $G$-diffeology might fail to be holonomy-like.

**Lemma 5.4.** The path holonomy diffeology on $R^G(\mathcal{B})$ is holonomy-like (see Def. 2.12).

*Proof.* Let $\chi\colon \mathcal{O} \to R^G(\mathcal{B})$ be a plot, and $e \subset \mathcal{O}$ a submanifold such that $\chi|_e\colon e \to \chi(e)$ is a diffeomorphisms onto an open subset of $1_M \subset R^G(\mathcal{B})$. We can lift $\chi$ to a plot $\tau$ for $H^G(\mathcal{B})$ (i.e. $\chi = \Phi \circ \tau$) after shrinking $\mathcal{O}$ if necessary, by the definition of path holonomy diffeology. Notice that $\tau(e)$ is a bisection of $H^G(\mathcal{B})$ that carries the identity of $G$. Thus $\tau(e)$ is contained in the set of units of $H^G(\mathcal{B})$, by the definition of holonomy groupoid (see the proof of Thm. 4.19 for more details).

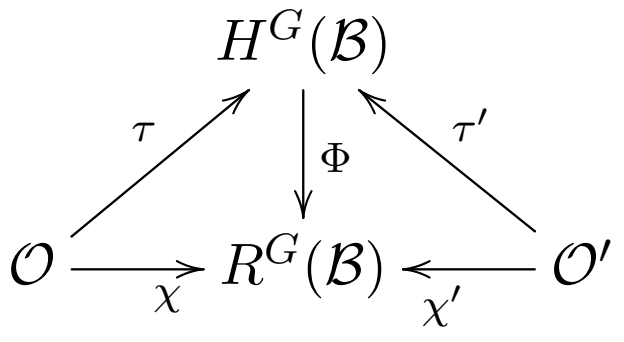

Repeat this procedure for another such plot $\chi'\colon \mathcal{O}' \to R^G(\mathcal{B})$ satisfying $\chi(e) = \chi'(e')$. Since $\tau|_e\colon e \to \tau(e)$ and $\tau'|_{e'}\colon e' \to \tau'(e')$ are diffeomorphisms onto the same open subset of $1_M \subset H^G(\mathcal{B})$, the fact that $H^G(\mathcal{B})$ is holonomy-like (Prop. 2.13) implies that there exists $k\colon \mathcal{O}' \to \mathcal{O}$ with $k(e') = e$ and $\tau \circ k = \tau'$, after shrinking the domain if necessary. Composing both sides with $\Phi$ we find $\chi \circ k = \chi'$. $\square$

**Proposition 5.5.** Endowed with the path holonomy diffeology, the graph $R^G(\mathcal{B})$ differentiates to $\mathcal{B}$ (Def. 3.3).

*Proof.* We will show that $\mathcal{B}$ equals

$$\left\{ \frac{d}{d\lambda}|_{\lambda=0}\mathbf{c}_\lambda\colon\ \{\mathbf{c}_\lambda\}_{\lambda\in I} \text{ smooth family of global bisections for } R^G(\mathcal{B}) \text{ s.t. } \mathbf{c}_0 = Id_M \right\} \cap \Gamma_c(AG).$$

Using this, it is clear that $R^G(\mathcal{B})$ satisfies Assumption 3.1, since $H^G(\mathcal{B})$ satisfies it and $\Phi$ is a smooth groupoid morphism. Together with Lemma 5.4, this implies the proposition.

"$\subset$:" Any $\boldsymbol{\alpha} \in \mathcal{B}$ is of the form $\frac{d}{d\lambda}|_{\lambda=0}(\Phi \circ \mathbf{b}_\lambda)$ for a smooth family of global bisections $\{\mathbf{b}_\lambda\}$ of $H^G(\mathcal{B})$, since the holonomy groupoid differentiates to $\mathcal{B}$ by Thm. A.7. Since $\mathbf{c}_\lambda := \Phi \circ \mathbf{b}_\lambda$ is a *smooth* family of global bisections of $R^G(\mathcal{B})$, this inclusion is proven.

"$\supset$:" let $\{\mathbf{c}_\lambda\}_{\lambda\in I}$ be smooth family of global bisections for $R^G(\mathcal{B})$ with $\mathbf{c}_0 = Id_M$, so that $\boldsymbol{\alpha} := \frac{d}{d\lambda}|_{\lambda=0}\mathbf{c}_\lambda$ lies in $\Gamma_c(AG)$. We have to show that $\boldsymbol{\alpha} \in \mathcal{B}$. To do this we use the fact that plots for $R^G(\mathcal{B})$ locally lift to plots for $H^G(\mathcal{B})$, by the definition of quotient diffeology. Take a finite open cover $\{U^i\}$ of $Supp(\boldsymbol{\alpha})$ so that the restriction of $\{\mathbf{c}_\lambda\}$ to $U^i$ lifts for every $i$, i.e. there is a smooth family of bisections $\{\mathbf{b}^i_\lambda\}$ for $H^G(\mathcal{B})$ defined on $U^i$ such that $\Phi \circ \mathbf{b}^i_\lambda = \mathbf{c}_\lambda|_{U^i}$ for all $\lambda$. Notice that $\mathbf{b}^i_0$ is contained in $1_M$, by the definition of holonomy groupoid. Take a partition of unity $\rho_i$ subordinate to this cover. Since $\rho_i\boldsymbol{\alpha} = \frac{d}{d\lambda}|_{\lambda=0}\mathbf{c}_{\rho_i\lambda}$ by the chain rule (exactly as in eq. (A.1)), it suffices to show that for all $i$

$$\frac{d}{d\lambda}|_{\lambda=0}\mathbf{c}_{\rho_i\lambda} \in \mathcal{B}. \tag{5.1}$$

Notice that $\mathbf{b}^i_{\rho_i\lambda}$ is a smooth family of bisections of $H^G(\mathcal{B})$. Since $\Phi \circ \mathbf{b}^i_{\rho_i\lambda} = \mathbf{c}_{\rho_i\lambda}|_{U^i}$ and the holonomy groupoid differentiates to $\mathcal{B}$, we deduce that eq. (5.1) holds. $\square$

By the following proposition, the graph $R^G(\mathcal{B})$ is generally fails to be an integral of $\mathcal{B}$. Indeed, it generally does not admit a generating set of plots consisting of open maps[13].

**Proposition 5.6.** The graph $R^G(\mathcal{B})$, endowed with the path holonomy diffeology, in general does *not* satisfy the openness condition b) in Def. 4.1. In particular, the graph $R^G(\mathcal{B})$ is generally not an integral.

We will prove this proposition as a consequence of the fact that the composition of an open map with the surjective map $\Phi : H^G(\mathcal{B}) \to R^G(\mathcal{B})$, where the latter has to quotient topology, is not an open map in general. Recall that the quotient topology agrees with the D-topology of the quotient diffeology, see Ex. 1.15 d).

*Proof.* We prove the proposition by providing a counter-example. Let $\mathcal{B} = \mathcal{F}$ be the singular foliation induced by the action of $K = S^1$ on $M = \mathbb{R}^2$ by rotations. We show that the graph $R(\mathcal{F}) \subset M \times M$ does not have any plot which is an open map with the point $(0,0) \in R(\mathcal{F})$ in its image.

Let $\chi \colon \mathcal{O} \to R(\mathcal{F})$ be a plot in the path holonomy diffeology and $x \in \mathcal{O}$ a point mapping to $(0,0)$. We want to show that $\chi$ is not an open map. There is a (open) neighborhood $\mathcal{O}^{(1)}$ of $x$ and a plot $\widetilde{\chi}$ for the holonomy groupoid making this diagram commute, by definition of path holonomy diffeology:

$$\begin{array}{ccc} \mathcal{O}^{(1)} & \xrightarrow{\widetilde{\chi}} & H(\mathcal{F}) \\ & \searrow^{\chi} & \downarrow \Phi \\ & & R(\mathcal{F}) \end{array} \tag{5.2}$$

Notice that the holonomy groupoid $H(\mathcal{F})$ is the transformation Lie groupoid $K \ltimes M$, and the map $\Phi$ is the target-source map $(k,v) \mapsto (kv, v)$.

Take a connected open subset $K^{(1)} \subsetneq K$ such that $K^{(1)} \times M$ contains $\widetilde{\chi}(x)$. Define $\mathcal{O}^{(2)}$ to be the preimage of $K^{(1)} \times M$ under $\widetilde{\chi}|_{\mathcal{O}^{(1)}}$. It is an open neighborhood of $x$ in $\mathcal{O}$, since $\widetilde{\chi}$ is continuous. The following claim immediately implies that $\chi$ is not an open map, as desired.

Claim: $\chi(\mathcal{O}^{(2)})$ is not an open subset of $R(\mathcal{F})$.

*Proof of claim:* The claim is equivalent to the following, since $R(\mathcal{F})$ is endowed with the quotient topology: $\Phi^{-1}(\Phi(\widetilde{\chi}(\mathcal{O}^{(2)}))$ is not an open subset of $H(\mathcal{F})$. Since by construction $\widetilde{\chi}(\mathcal{O}^{(2)}) \subset K^{(1)} \times M$, we have

$$\Phi^{-1}\left(\Phi(\widetilde{\chi}(\mathcal{O}^{(2)})\right) \subset \Phi^{-1}\left(\Phi(K^{(1)} \times M)\right). \tag{5.3}$$

The right hand side equals $(K^{(1)} \times M \backslash \{0\}) \coprod K \times \{0\}$, since the $K$-action is free on $M\backslash\{0\}$ and fixes $\{0\}$. Notice that the left hand side contains $K \times \{0\}$. Pick $k \in K \backslash K^{(1)}$. Then $(k,0)$ lies in the left hand side, but due to the sentence following (5.3) it does not admit any open neighborhood in $H(\mathcal{F})$ lying in the right hand side. △

□

We now highlight properties of the diffeological groupoid $H^G(\mathcal{B})$ that are not shared by $R^G(\mathcal{B})$.

[13] The graph $R^G(\mathcal{B})$ does satisfy the weaker property of being source-submersive, by the definition of path-holonomy diffeology, Cor. 2.25 and Prop. 2.17.

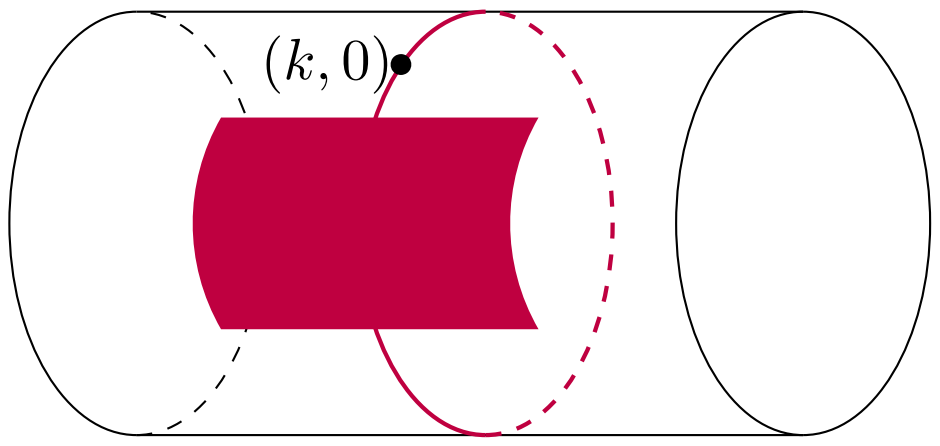

Figure 1: The point $(k,0)$ lies on the subset $(K^{(1)} \times M\backslash\{0\}) \coprod K \times \{0\}$ (in purple color) of $K \ltimes M$. Here $M$ is represented as an interval, for easiness of drawing.

- Integrals: We saw that $H^G(\mathcal{B})$ is an integral of $\mathcal{B}$, while $R^G(\mathcal{B})$ is not in general (Thm. 4.19 and Prop. 5.6).
- Smoothness: In the regular case, i.e. when $\mathcal{B}$ comes from a wide Lie subalgebroid $B \subset AG$, the holonomy groupoid $H^G(\mathcal{B})$ is smooth, while in general $R^G(\mathcal{B})$ fails to be a Lie groupoid integrating $B$.
- Dimension of source-fibers: The dimension of the source-fibers of $H(\mathcal{B})$ is upper semicontinuous, i.e. locally it can not increase. In contrast to this, the dimension of the source-fibers of $R(\mathcal{B})$ is lower semicontinuous, as is apparent in the case of singular foliations.
- Holonomy transformations: It was shown in [7] that the holonomy groupoid $H(\mathcal{F})$ of a singular foliation has a canonical injective morphism into the groupoid of holonomy transformations. The latter are the germs of diffeomorphisms between (fixed) slices transverse to $\mathcal{F}$, modulo a certain equivalence relation. The injectivity implies that this morphism does not descend to any (non-trivial) quotient of $H(\mathcal{F})$, in particular not to $R(\mathcal{F})$.

  These statements extend to any singular subalgebroid of a Lie algebroid $AG$, as shown in [30]. (The slices are then chosen as follows: the slice through a point $x \in M$ is contained in the source-fiber $\mathbf{s}^{-1}(x) \subset G$ and is transverse to the evaluation of $\mathcal{B}$ at $x$.)

As seen in Thm. 4.19, the holonomy groupoid is always an integral. Adapting the proof of Prop. 5.6, we present an example in which all its non-trivial quotients fail to be an integral.

**Corollary 5.7.** Let $\mathcal{F}$ be the singular foliation induced by the action of $S^1$ on $\mathbb{R}^2$ by rotations. Then no non-trivial quotient of $H(\mathcal{F})$ to a groupoid over $\mathbb{R}^2$, with the quotient diffeology, is an integral of $\mathcal{F}$.

*Proof.* The holonomy groupoid $H(\mathcal{F})$ is the transformation Lie groupoid $K \ltimes M$, where $K = S^1 = U(1)$ and $M = \mathbb{R}^2$. We consider quotients of $H(\mathcal{F})$ to groupoids over the same manifold of units $M$. Since the $K$-action is free on $M\backslash\{0\}$, the quotient can only affect the restriction of $H(\mathcal{F})$ to the origin $0 \in M$, which is the abelian group $K$. For the sake of presentation, we consider a specific quotient of $H(\mathcal{F})$, but the argument holds for any quotient as above. Let $H'$ be the quotient of $H(\mathcal{F})$ by the equivalence relation given by the antipodal $\mathbb{Z}_2$-action on the fiber over the origin, i.e. $k \times \{0\} \sim (-k) \times \{0\}$ for all $k \in K$ (over $M\backslash\{0\}$ the equivalence relation is trivial). Then $H'$ is again a diffeological groupoid, when endowed with the quotient diffeology via the quotient map $\pi \colon H(\mathcal{F}) \to H'$. Notice that even though the kernel of $\pi$ is discrete, neighbourhoods of the units in $H(\mathcal{F})$ and $H'$ look quite different: for instance, $e^{i(-\epsilon,\epsilon)} \times M$ is an open subset of $H(\mathcal{F})$, but its image under $\pi$ is not open in $H'$ with the quotient topology, because $\pi^{-1}(\pi(e^{i(-\epsilon,\epsilon)} \times M))$ also contains the "arc" $-e^{i(-\epsilon,\epsilon)} \times \{0\}$.

We show that the openness condition b) in Def. 4.1 is not satisfied by $H'$. Indeed, we argue that $H'$ does not have any plot which is an open map with the point $\pi(0,1) = \pi(0,-1)$ in its image. To this aim, we just have to modify the proof of Prop. 5.6, as follows. Take the open subset $K^{(1)} \subsetneq K$ there to be small enough. The right hand side of (5.3) (writing $\pi$ instead of $\Phi$) reads $(K^{(1)} \times M) \coprod (-K^{(1)} \times \{0\})$. Notice that the left hand side of (5.3) contains $k := -\tilde{\chi}(x) \in -K^{(1)}$ Then $(k, 0)$ lies in the left hand side, but it does not admit any open neighborhood in $H(\mathcal{F})$ lying in the right hand side. □

## 5.2 The graph as image of arbitrary integrals

Let $\mathcal{B}$ be a singular subalgebroid of a Lie algebroid $AG$. We show that for all integrals of $\mathcal{B}$, the image of the natural map to $G$ is always the same, namely the graph $R^G(\mathcal{B})$. We provide a condition under which the quotient diffeology on the image coincides with the path-holonomy diffeology introduced in §5.1.2. We do now know if the condition is always satisfied or not; this question is related to the question of whether the holonomy groupoid is a minimal integral.

**Proposition 5.8.** Let $(H, \Psi)$ be an integral of $\mathcal{B}$ over $G$. Then $\Psi(H) = R^G(\mathcal{B})$.

In the proof of this proposition, the inclusion "$\supset$" makes essential use of the coarsely injective condition c) in Def. 4.1, and the inclusion "$\subset$" of the openness condition b) there.

*Proof.* "$\supset$": Since $(H, \Psi)$ is an integral, by Cor. B.2 it satisfies Assumption B.3. Notice that Assumption B.3 requires for every $x \in M$ to choose a neighborhood, which we denote $V_x$. Endow $R^G(\mathcal{B})$ with the quotient topology induced by $\Phi\colon H^G(\mathcal{B}) \to R^G(\mathcal{B})$.

Claim: *There is a neighborhood $U$ of $1_M$ in $R^G(\mathcal{B})$ such that for every $g \in U$ there is $\boldsymbol{\alpha} \in \mathcal{B}$ with $g \in \exp_{1_M}(\overrightarrow{\boldsymbol{\alpha}})$ and $Supp(\boldsymbol{\alpha}) \subset V_{\mathbf{s}(g)}$.*

*Proof of claim:* Let $x \in M$, and denote by $L \subset M$ the leaf of $\mathcal{B}$ through $x$. Assume first that $L$ is embedded. In view of Remark 5.2, consider the transitive Lie algebroid $B_L$ over $L$ characterized by $\Gamma_c(B_L) \cong \{\boldsymbol{\alpha}|_L : \boldsymbol{\alpha} \in \mathcal{B}\}$. It integrates the (transitive) Lie groupoid $R^G(\mathcal{B})_L$. It is a general fact about Lie groupoids that there is a neighborhood of $1_L$ in $R^G(\mathcal{B})_L$ such that every point $g$ there lies on $\exp_{1_L}(\overrightarrow{a})$ for some $a \in \Gamma_c(B_L)$ [17, §4, page 16]. We will use the following stronger statement [32, Lemma 5.2]: for every neighborhood $V_x^L$ of $x$ in $L$ there is an open neighborhood $U_x$ of $1_x$ in $\mathbf{s}^{-1}_{R^G(\mathcal{B})_L}(x)$ such that any point $q \in U_x$ lies on $\exp_{1_L}(\overrightarrow{a})$ for some $a \in \Gamma_c(B_L)$ with $supp(a) \subset V_x^L$. We apply this to $V_x^L = V_x \cap L$. There exists $\boldsymbol{\alpha} \in \mathcal{B}$ with $\boldsymbol{\alpha}|_L = a$, which can be chosen so that $Supp(\boldsymbol{\alpha}) \subset V_x$. By construction, $g \in \exp_{1_M}(\overrightarrow{\boldsymbol{\alpha}})$.

If the leaf $L$ is not embedded in $M$, we can apply the argument above to a neighborhood $L_x \subset L$ which is an embedded submanifold of $M$. To finish the argument, take $U$ to be a neighborhood of $1_M$ contained in $\cup_{x \in M} U_x \subset R^G(\mathcal{B})$. △

Now let $g \in U \subset R^G(\mathcal{B})$ be as in the claim. By Assumption B.3 there are $N \in \mathbb{N}$ and a 1-parameter family $\{\mathbf{b}_\lambda\}_{\lambda \in [0,\frac{1}{N}]}$ of global bisections of $H$ (with $\mathbf{b}_0 = 1_M$) mapping to $\{\exp_{1_M}(\lambda\overrightarrow{\boldsymbol{\alpha}})\}_{\lambda \in [0,\frac{1}{N}]}$ under $\Psi$. In particular, for $x := \mathbf{s}(g)$, we have $\Psi\left(\mathbf{b}_{\frac{1}{N}}(x)\right) = \exp_x\left(\frac{1}{N}\overrightarrow{\boldsymbol{\alpha}}\right)$, and taking the $N$-th power we see that $\exp_x(\boldsymbol{\alpha}) = g \in \Psi(H)$. Since the groupoid $R^G(\mathcal{B})$ is source-connected, the desired inclusion $\Psi(H) \supset R^G(\mathcal{B})$ follows.

"⊂" We proceed as in the proof of Prop. 4.7. By the openness condition b) in Def. 4.1 we can apply Prop. 2.17 ii). Hence there is a neighborhood $\breve{H}$ of the identity $1_M$ in $H$ such that for all $h \in \breve{H}$ there is a 1-parameter family[14] of global bisections $\{\mathbf{b}_\lambda\}_{\lambda\in[0,1]}$ with $\mathbf{b}_0 = 1_M$ and $h \in \mathbf{b}_1$. (This family can be chosen with compact support). Let $h \in \breve{H}$ and $\{\mathbf{b}_\lambda\}_{\lambda\in[0,1]}$ be as above. For every fixed $\lambda_0 \in [0,1]$ consider $\{\mathbf{b}_{\mu+\lambda_0} * \mathbf{b}_{\lambda_0}^{-1}\}_{\mu\in I}$ where $I$ is an open interval containing zero. It is a 1-parameter family of bisections of $H$. Since $H$ differentiates to $\mathcal{B}$, we therefore have

$$\frac{d}{d\mu}|_{\mu=0}\left(\Psi \circ (\mathbf{b}_{\mu+\lambda_0} * \mathbf{b}_{\lambda_0}^{-1})\right) \in \mathcal{B}.$$

Since $\Psi$ is a groupoid morphism, the left hand side equals

$$\frac{d}{d\mu}|_{\mu=0}R_{\Psi\circ\mathbf{b}_{\lambda_0}^{-1}}(\Psi \circ \mathbf{b}_{\mu+\lambda_0}) = (R_{\Psi\circ\mathbf{b}_{\lambda_0}^{-1}})_*\frac{d}{d\lambda}|_{\lambda_0}(\Psi \circ \mathbf{b}_\lambda)$$

where $R$ denotes right-translation on $G$. Since $\overrightarrow{\mathcal{B}}$ is right-invariant (see §1.1), from this we conclude that $\frac{d}{d\lambda}|_{\lambda_0}(\Psi \circ \mathbf{b}_\lambda)$ is the restriction of an element of $\overrightarrow{\mathcal{B}}$ to $\Psi \circ \mathbf{b}_{\lambda_0}$. In particular, for $x := \mathbf{s}_H(h)$, the curve $[0,1] \ni \lambda \mapsto (\Psi \circ \mathbf{b}_\lambda)(x)$ is tangent to $\overrightarrow{\mathcal{B}}$, so its endpoint $\Psi(h)$ lies in $R^G(\mathcal{B})$. □

We now address diffeological structures. Thanks to Prop. 5.8, we know that $\Psi(H)$ agrees $R^G(\mathcal{B}) = \Phi(H^G(\mathcal{B}))$ as a set, so the assumption of the following statement is reasonable.

**Proposition 5.9.** Let $(H, \Psi)$ be an integral of $\mathcal{B}$ over $G$. Assume that there is a subduction $\pi$ to the holonomy groupoid such that the following diagram commutes:

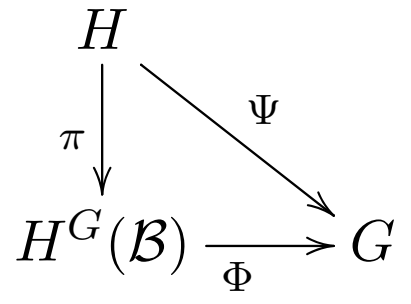

Then the quotient diffeology on $\Psi(H)$ (induced by $\Psi\colon H \to \Psi(H)$) agrees with the path-holonomy diffeology.

We recall from §1.4 what it means that $\pi\colon H \to H^G(\mathcal{B})$ is a subduction: it means that $\pi$ is a surjective map and that the quotient diffeology induced by it coincides with the diffeology on $H^G(\mathcal{B})$.

*Proof.* By definition, the path-holonomy diffeology is the quotient diffeology obtained from the surjective map $H^G(\mathcal{B}) \to \Phi(H^G(\mathcal{B})) = R^G(\mathcal{B})$. The quotient diffeology is well-behaved under composition [16, §1.45], so the commutativity of the above diagram implies the statement. □

# A On differentiation

In §3 we formulated a definition of differentiation, and argued that it is a reasonable one based on Thm. 3.2. In this appendix we provide proofs for Thm. 3.2 and for the two examples displayed in §3. Further, we make remarks on the functoriality of differentiation.

[14] We clarify what we mean by 1-parameter family parametrized by the closed interval $[0,1]$: by definition, this means that there is an open interval $I$ containing $[0,1]$ and a (smooth) 1-parameter family defined in $I$ extending $\{\mathbf{b}_\lambda\}_{\lambda\in[0,1]}$.

## A.1 The proof of theorem 3.2

Let $H \rightrightarrows M$ be a diffeological groupoid which is holonomy-like, $G \rightrightarrows M$ a Lie groupoid and $\Psi : H \to G$ a smooth morphism of diffeological groupoids covering $Id_M$ and satisfying Assumption 3.1. Recall that the subset $\mathcal{S}$ of $\Gamma_c(AG)$ was defined in §3.

The proof of theorem 3.2 – which states that $\mathcal{S}$ is a singular subalgebroid – is provided by the following three lemmas.

**Lemma A.1.** The set $\mathcal{S} \subseteq \Gamma_c(AG)$ is a $C^\infty(M)$-module.

*Proof.* Given $f \in C^\infty(M)$ and an element $\frac{d}{d\lambda}|_{\lambda=0}(\Psi \circ \mathbf{b}_\lambda)$ of $\mathcal{S}$, we show that their product lies in $\mathcal{S}$. Since elements of $\mathcal{S}$ are compactly supported, we may assume that $f$ is compactly supported. Consider the family of global bisections defined by $\mathbf{b}^f_\lambda(x) = \mathbf{b}_{f(x)\lambda}$. (Since $f$ is a bounded function, $\lambda$ is defined in an open subinterval of $I$ containing zero, and passing to a smaller subinterval if necessary we ensure that $\mathbf{b}^f_\lambda$ is a bisection.) The chain rule implies

$$\frac{d}{d\lambda}|_{\lambda=0}(\Psi \circ \mathbf{b}^f_\lambda) = f\frac{d}{d\lambda}|_{\lambda=0}(\Psi \circ \mathbf{b}_\lambda), \tag{A.1}$$

and the l.h.s. clearly lies in $\mathcal{S}$.

Now fix two elements $\frac{d}{d\lambda}|_{\lambda=0}(\Psi \circ \mathbf{b}_\lambda)$ and $\frac{d}{d\lambda}|_{\lambda=0}(\Psi \circ \mathbf{c}_\lambda)$ of $\mathcal{S}$, where $\mathbf{b}_\lambda$ and $\mathbf{c}_\lambda$ are 1-parameter families with $\mathbf{b}_0 = \mathbf{c}_0 = Id_M$. We show that their sum lies in $\mathcal{S}$, and more concretely that it comes from the 1-parameter family $\mathbf{b}_\lambda * \mathbf{c}_\lambda$ of global bisections of $H$. Indeed, using the product of global bisections introduced in §2.2, we have

$$\frac{d}{d\lambda}|_{\lambda=0}(\Psi \circ (\mathbf{b}_\lambda * \mathbf{c}_\lambda)) = \frac{d}{d\lambda}|_{\lambda=0}((\Psi \circ \mathbf{b}_\lambda) * (\Psi \circ \mathbf{c}_\lambda)) = \frac{d}{d\lambda}|_{\lambda=0}((\Psi \circ \mathbf{b}_\lambda) + \frac{d}{d\lambda}|_{\lambda=0}(\Psi \circ \mathbf{c}_\lambda),$$

where in the first equality we used that $\Psi$ is a groupoid morphism and in the second that $\mathbf{b}_0 = \mathbf{c}_0$ is the identity bisection. ☐

The assumption that the $H$ is holonomy-like is used to prove that $\mathcal{S}$ is locally finitely generated.

**Lemma A.2.** The module $\mathcal{S} \subseteq \Gamma_c(AG)$ is locally finitely generated.

*Proof.* For every $x \in 1_M$, by Lemma 2.11 we know that there is a neighbourhood $V$ of $x$ in $M$, a plot $\chi$ in a generating set $\mathcal{X}$ as in Def. 2.12, and a submanifold $e$ of $\mathcal{O}_\chi$, such that

- $\chi|_e \colon e \to 1(V) \subset 1_M$ is a diffeomorphism,
- $\mathbf{s}_H \circ \chi$ and $\mathbf{t}_H \circ \chi$ are a submersion at points of $e$.

For later use, we remark that $e$ is realized as the image of a smooth map $\tilde{1} \colon V \to \mathcal{O}_\chi$.

$$\begin{array}{ccc} e \subset \mathcal{O}_\chi & \xrightarrow{\chi} & H \\ & \nwarrow^{\tilde{1}} & \uparrow{\scriptstyle 1|_V} \\ & & V \end{array} \tag{A.2}$$

The above implies that

$$Vert|_e := Ker(\mathbf{s}_H \circ \chi)_*|_e$$

is a vector bundle over $e$.

<u>Claim:</u>

$$\mathcal{S}|_V = (\Psi \circ \chi)_* \Gamma(Vert|_e).$$

Since $\Gamma(Vert|_e)$ is finitely generated as a $C^\infty(e)$-module, the claim implies that the same holds for $\mathcal{S}|_V$, as we wanted to show. Thus we are left with proving the claim.

*Proof of claim:*

"$\subset$" Every element of $\mathcal{S}|_V$ is of the form $\frac{d}{d\lambda}|_0(\Psi \circ \mathbf{b}_\lambda)$ where $\{\mathbf{b}_\lambda\}_{\lambda \in I}$ is a 1-parameter family of bisections of $H$ defined on $V$, with $\mathbf{b}_0 = 1|_V$. By Def. 2.6, $\mathbf{b} : I \times V \to H$ is smooth.

Fix $y \in V$. There is an open neighborhood $V_y$ in $V$ such that the map $\mathbf{b} : I \times V \to \overset{\circ}{H}$ is well-defined, where $\overset{\circ}{H}$ is as in Def. 2.12, shrinking $I$ to a smaller open interval about zero if necessary. Shrinking $V_y$ if necessary, $\mathbf{b}|_{I \times V_y}$ is a plot in $\mathcal{P}(\overset{\circ}{H})$, by item c) in §1.4 and by Prop. 1.12, and there exists a plot $\chi' : \mathcal{O}_{\chi'} \to \overset{\circ}{H}$ in the generating set $\mathcal{X}$ (as in Def. 2.12) and a smooth map $h' : I \times V_y \to \mathcal{O}_{\chi'}$ such that $\mathbf{b}|_{I \times V_y} = \chi' \circ h'$. For the sake of readability, let us temporarily assume that $V_y = V$.

$$\begin{array}{ccc} \mathcal{O}_{\chi'} & \xrightarrow{\chi'} & H \\ & \nwarrow^{h'} & \uparrow \mathbf{b} \\ & & I \times V \end{array} \tag{A.3}$$

We just "lifted" $\mathbf{b}$ to $\mathcal{O}_{\chi'}$, and now we argue that $\mathbf{b}$ can be lifted to the fixed plot $\mathcal{O}_\chi$ introduced at the beginning of the proof, in such a way that $\{0\} \times V$ is mapped to the submanifold $e$ fixed there. Notice that $\chi'$ maps the submanifold $h'_0(V)$ diffeomorphically onto $\mathbf{b}_0(V) = 1(V)$, so since $H$ is holonomy-like (Def. 2.12) there exists a smooth map $k \colon \mathcal{O}_{\chi'} \to \mathcal{O}_\chi$ such that $\chi \circ k = \chi'$ mapping $h'_0(V)$ to $e$, shrinking $\mathcal{O}_{\chi'}$ to a smaller neighborhood of $h'_0(y)$ if necessary. Therefore

$$h := k \circ h' \colon I \times V \to \mathcal{O}_\chi$$

satisfies $\chi \circ h = \chi \circ k \circ h' = \chi' \circ h' = \mathbf{b}$ (i.e., the diagram below commutes) and $h_0(V) = k(h'_0(V)) = e$.

$$\begin{array}{ccc} \mathcal{O}_{\chi} & \xrightarrow{\chi} & H \\ & \nwarrow^{h} & \uparrow \mathbf{b} \\ & & I \times V \end{array} \tag{A.4}$$

Even more, $h_0$ coincides with $\tilde{1} \colon V \to e$, by the commutativity of diagrams (A.2) and (A.4) together with the facts that $1|_V$ and $\chi|_e$ are diffeomorphisms onto their images.

For every $\lambda \in I$, $h_\lambda$ is a section of $\mathbf{s}_H \circ \chi$, since $\mathbf{b}_\lambda = \chi \circ h_\lambda$ is a section of $\mathbf{s}_H$. Consider $\Psi \circ \mathbf{b}_\lambda = (\Psi \circ \chi) \circ h_\lambda$. Taking derivatives (notice that $\Psi \circ \chi \colon \mathcal{O}_\chi \to G$ is smooth) we obtain

$$\frac{d}{d\lambda}|_0(\Psi \circ \mathbf{b}_\lambda) = (\Psi \circ \chi)_* \frac{d}{d\lambda}|_{\lambda=0}(h_\lambda). \tag{A.5}$$

Notice that $\frac{d}{d\lambda}|_{\lambda=0}(h_\lambda)$ is a section of $Vert|_e$, proving the desired inclusion under the assumption $V_y = V$.

For the general case, we argue as follows. There are open subsets $V^i$ covering $V$ and maps $h^i\colon I \times V^i \to \mathcal{O}_\chi$ so that $(h^i)_0\colon V^i \to \mathcal{O}_\chi$ equals $\tilde{1}|_{V^i}$. Our reasoning above shows that eq. (A.5) holds on $V^i$, replacing $h$ by $h^i$. Let $\{\varphi^i\}$ be a partition of unity subordinate to the open cover $\{\tilde{1}(V^i)\}$ of $e$. Then eq. (A.5) holds, on the whole of $V$, replacing the section $\frac{d}{d\lambda}|_{\lambda=0}(h_\lambda)$ of $Vert|_e$ with the section $\sum_i(\tilde{1}^*\varphi^i)\frac{d}{d\lambda}|_{\lambda=0}(h^i_\lambda)$.

" $\supset$ Any element of $Vert|_e$ can be written as $\frac{d}{d\lambda}|_{\lambda=0}(h_\lambda)$ for a smooth family of sections $h_\lambda\colon V \to \mathcal{O}_\chi$ of $\mathbf{s}_H \circ \chi$ through $e$. (This uses the fact that since $\mathbf{s}_H \circ \chi$ is a submersion at points of $e$, it is a submersion also in a tubular neighbourhood of it.) Further $\chi \circ h_\lambda$ is a family of bisections of $H$ which is smooth by construction and goes through $1_V$. Hence $\frac{d}{d\lambda}|_{\lambda=0}(\Psi \circ \chi \circ h_\lambda)$ is an element of $\mathcal{S}|_V$. △

□

Assumption 3.1 is used to prove that $\mathcal{S}$ is involutive. We will make use of the following fact about 1-parameter *groups* of global bisections of Lie groupoids. Let $\{\sigma_\lambda : M \to G\}_{\lambda\in I}$ be such a 1-parameter group. For any fixed $\lambda \in I$, denote by

$$L_{\sigma_\lambda}\colon G \to G, g \mapsto \sigma_\lambda(\mathbf{t}_G(g)) \cdot g$$

the left translation by $\sigma_\lambda$. Then $L_{\sigma_{\lambda+\nu}} = L_{\sigma_\lambda * \sigma_\nu} = L_{\sigma_\lambda} \circ L_{\sigma_\nu}$, therefore $\{L_{\sigma_\lambda}\}_{\lambda\in I}$ is the flow of a vector field on $G$. By differentiation it is easily seen that this vector field is the right-invariant vector field $\overrightarrow{\frac{d}{d\lambda}|_0\sigma_\lambda}$.

**Lemma A.3.** The module $\mathcal{S} \subseteq \Gamma_c(AG)$ is involutive.

*Proof.* Let $\boldsymbol{\alpha}_1, \boldsymbol{\alpha}_2 \in \mathcal{S}$. Let $\{V_i\}$ be a cover of $M$ by open subsets on which the inclusion in Assumption 3.1 holds; we may assume that finitely many of these open sets cover $Supp(\boldsymbol{\alpha}_1) \cap Supp(\boldsymbol{\alpha}_2)$. Let $\{\varphi_i\}$ a partition of unity subordinate to $\{V_i\}$. Then $[\boldsymbol{\alpha}_1, \boldsymbol{\alpha}_2] = \sum_{i,j}[\varphi_i\boldsymbol{\alpha}_1, \varphi_j\boldsymbol{\alpha}_2]$. Thanks to Lemma A.1, we only need to show that each summand lies in $\mathcal{S}$.

Fix indices $i, j$. Thanks to Lemma A.1 and Assumption 3.1 we can write $\varphi_i\boldsymbol{\alpha}_1 = \frac{d}{d\lambda}|_0(\Psi \circ \mathbf{b}_\lambda)$ and $\varphi_j\boldsymbol{\alpha}_2 = \frac{d}{d\lambda}|_0(\Psi \circ \mathbf{b}'_\lambda)$, where $\{\mathbf{b}_\lambda\}_{\lambda\in I}, \{\mathbf{b}'_\lambda\}_{\lambda\in I}$ are 1-parameter *groups* of global bisections of $H$. Their Lie bracket is the restriction to $M$ of the Lie bracket of the associated right-invariant vector fields on $G$

$$\left[\overrightarrow{\frac{d}{d\lambda}|_0(\Psi \circ \mathbf{b}_\lambda)}\, , \,\overrightarrow{\frac{d}{d\lambda}|_0(\Psi \circ \mathbf{b}'_\lambda)}\right], \tag{A.6}$$

whose flows are given by the left translations $\{L_{\Psi\circ\mathbf{b}_\lambda}\}_{\lambda\in I}$ and $\{L_{\Psi\circ\mathbf{b}'_\lambda}\}_{\lambda\in I}$ respectively. Using the characterization of the Lie bracket of two vector fields on any manifold in terms of their flows, we find that the value of (A.6) at any $g \in G$ is

$$\frac{d}{d\lambda}|_0 \left( L_{\Psi\circ\mathbf{b}'_{-\sqrt{\lambda}}} \circ L_{\Psi\circ\mathbf{b}_{-\sqrt{\lambda}}} \circ L_{\Psi\circ\mathbf{b}'_{\sqrt{\lambda}}} \circ L_{\Psi\circ\mathbf{b}_{\sqrt{\lambda}}} \right)(g) =$$
$$\frac{d}{d\lambda}|_0 \left( L_{\Psi\circ(\mathbf{b}'_{-\sqrt{\lambda}} * \mathbf{b}_{-\sqrt{\lambda}} * \mathbf{b}'_{\sqrt{\lambda}} * \mathbf{b}_{\sqrt{\lambda}})} \right)(g).$$

Evaluating at all points of $M$ this yields an element of $\mathcal{S}$, because $\mathbf{b}'_{-\sqrt{\lambda}} * \mathbf{b}_{-\sqrt{\lambda}} * \mathbf{b}'_{\sqrt{\lambda}} * \mathbf{b}_{\sqrt{\lambda}}$ is a 1-parameter family of global bisections of $H$. Above we use the fact that $\Psi$, being a groupoid morphism covering $Id_M$, respects the multiplication $*$ of 1-parameter groups of global bisections Thus $[\varphi_i\boldsymbol{\alpha}_1, \varphi_j\boldsymbol{\alpha}_2] \in \mathcal{S}$, concluding the proof. □

Theorem 3.2 is thus proven, thanks to Lemmas A.1, A.2, A.3.

### A.1.1 About the module $\mathcal{S}$

We end this subsection arguing that $\mathcal{S}$ "contains many elements". First notice that any $C^\infty(M)$-module of compactly supported sections is generated by its elements supported on arbitrarily small open subsets, as one can prove using a partition of unity argument:

**Remark A.4.** Let $\mathcal{M} \subseteq \Gamma_c(AG)$ be a $C^\infty(M)$-module and $\{V_i\}$ be open cover of $M$. Then $\mathcal{M}$ is generated by its element whose support is contained in $V_i$ for some index $i$.

The following proposition provides elements of $\mathcal{S}$ supported on arbitrarily small open subsets. It shows that if a compactly supported section $\boldsymbol{\alpha}$ of $AG$ arises from a family of *not necessarily compactly supported* bisections of $H$, then $\boldsymbol{\alpha}$ lies in $\mathcal{S}$.

**Proposition A.5.** Let $V \subset M$ be an open subset, let $\boldsymbol{\alpha} = \frac{d}{d\lambda}|_{\lambda=0}(\Psi \circ \mathbf{b}_\lambda) \in \Gamma_c(AG|_V)$ where $\{\mathbf{b}_\lambda\}_{\lambda\in I}$ is family of bisections of $H$ defined on $V$ s.t. $\mathbf{b}_0 = Id_V$. Then one can choose the $\mathbf{b}_\lambda$'s so that their support is contained in a compact subset $K$ of $V$, shrinking $I$ if necessary. Therefore $\boldsymbol{\alpha}$, extended trivially to the rest of $M$, lies in $\mathcal{S}$.

*Proof.* Take $f \in C_c^\infty(V)$ which is identically equal to 1 on $supp(\boldsymbol{\alpha})$. Then

$$\boldsymbol{\alpha} = f\boldsymbol{\alpha} = \frac{d}{d\lambda}|_{\lambda=0}(\Psi \circ \mathbf{b}^f_\lambda)$$

by eq. (A.1). Notice that the bisections $\mathbf{b}^f_\lambda := \mathbf{b}_{f\lambda}$ have support contained in $K := supp(f)$ (meaning that they are the identity bisection outside of there). Extending these bisections of $H$ to global bisections, the last statement of the proposition follows. □

## A.2 Differentation and Lie groupoids

A first class of examples for differentiation (Def. 3.3) is provided by Lie groupoids:

**Proposition A.6.** Let $H$ be a source-connected Lie groupoid and $\Psi\colon H \to G$ a Lie groupoid morphism covering $Id_M$. Define $\mathcal{B} := \Psi_*(\Gamma_c(AH))$. Then $H$ differentiates to $\mathcal{B}$.

*Proof.* Every Lie groupoid is holonomy-like, by Remark 2.14. Since $H$ is a Lie groupoid, we have

$$\begin{aligned}\Gamma(AH) &= \left\{ \frac{d}{d\lambda}|_{\lambda=0}\mathbf{b}_\lambda \colon \ \{\mathbf{b}_\lambda\}_{\lambda\in I} \text{ 1-parameter } \textit{group} \text{ of global bisections of } H \right\} \\ &= \left\{ \frac{d}{d\lambda}|_{\lambda=0}\mathbf{b}_\lambda \colon \ \{\mathbf{b}_\lambda\}_{\lambda\in I} \text{ 1-parameter family of global bisections of } H \text{ with } \mathbf{b}_0 = Id_M \right\},\end{aligned}$$

where the first equality follows using arguments similar to those just before Lemma A.3. Hence $\mathcal{S}$ equal $\Psi_*(\Gamma_c(AH)) = \mathcal{B}$. Further, the first equality implies that Assumption 3.1 is satisfied. □

## A.3 Differentiation of the holonomy groupoid

A further class of examples for differentiation is given by holonomy groupoids. Let $\mathcal{B}$ be a singular subalgebroid of a Lie algebroid $A$ with source-connected integration $G$.

**Theorem A.7.** The holonomy groupoid $H^G(\mathcal{B})$ differentiates to $\mathcal{B}$ (in the sense of definition 3.3).

*Proof.* Since the holonomy-like property holds by Prop. 2.13, we only have to show that $\Phi : H^G(\mathcal{B}) \to G$ satisfies $\mathcal{S} = \mathcal{B}$ and Assumption 3.1.

We first prove $\mathcal{B} \subset \mathcal{S}$. Since $\mathcal{B}$ is a $C^\infty(M)$-submodule of $\Gamma_c(AG)$, by Remark A.4 it is generated by its sections which are supported on open neighborhoods chosen as in Proposition 2.9. Applying Proposition 2.9 and using the fact that $\mathcal{S}$ is a $C^\infty(M)$-submodule (by Lemma A.1) one obtains the inclusion $\mathcal{B} \subset \mathcal{S}$.

Now we prove $\mathcal{S} \subset \mathcal{B}$. To this end, we need some preparation. Let $x \in M$, and let $V$ be a sufficiently small neighborhood of $x$. Recall that:

- Given any smooth bisection $\mathbf{b}$ of $H^G(\mathcal{B})$ defined on $V$ and close enough (in the $C^0$-sense) to the identity bisection, there is a minimal path-holonomy bisubmersion $(U, \varphi, G)$ and a bisection $\mathbf{u}$ of $U$ such that $q_U \circ \mathbf{u} = \mathbf{b}$, by Rem. 2.4 a). Notice that $\varphi \circ \mathbf{u} = \Phi \circ \mathbf{b}$.
- Denote by $\boldsymbol{\alpha}_1, \ldots, \boldsymbol{\alpha}_n \in \mathcal{B}$ the local generators of $\mathcal{B}$ (inducing a basis of $\mathcal{B}/I_x\mathcal{B}$) used to construct the path-holonomy bisubmersion $U \subset \mathbb{R}^n \times M$. For any bisection $\mathbf{u} = (Id, f)$ of $U$ (where $f \colon M \to \mathbb{R}^n$), we have $(\varphi \circ \mathbf{u})(y) = \exp_{(0,y)}(\sum f_i(y)\overrightarrow{\boldsymbol{\alpha}_i})$, see item a) in Ex. 1.5.

Let $\{\mathbf{b}_\lambda\}$ be a smooth 1-parameter family of bisections of $H^G(\mathcal{B})$ defined on $V$, such that $\mathbf{b}_0 \subset Id_M$. Shrinking $V$ if necessary, we find a family of bisections $\mathbf{u}_\lambda$ of $U$ such that $q_U \circ \mathbf{u}_\lambda = \mathbf{b}_\lambda$, and one can always arrange[15] that $\mathbf{u}_0$ is the zero section of $U$.

Write $\mathbf{u}_\lambda = (Id, f_\lambda)$. Then for all $y \in V$ we have

$$\frac{d}{d\lambda}|_{\lambda=0}(\varphi \circ \mathbf{u}_\lambda)(y) = d_{(0,y)}\varphi(\frac{d}{d\lambda}|_{\lambda=0}\mathbf{u}_\lambda(y)) = \sum_i g_i(y)\overrightarrow{\boldsymbol{\alpha}_i}|_y, \quad \text{for} \quad g(y) = \frac{d}{d\lambda}|_{\lambda=0}(f_\lambda)(y),$$

[15] Indeed, since $q_U \circ \mathbf{u}_0 = \mathbf{b}_0$ is the identity bisection of $H^G(\mathcal{B})$, composing with $\Phi$ we see that $\varphi \circ \mathbf{u}_0$ is the identity section of $G$. Hence there is a (locally defined) morphism of bisubmersions $k \colon U \to U$ mapping $\mathbf{u}_0(V)$ to $(0, V)$, by [5, Cor. 3.3] (see the last three items in §1.2.1) and its proof. Since $q_U \circ k = q_U$, it follows that $k \circ \mathbf{u}$ is also a family of bisections of $U$ lifting $\mathbf{b}$, with the additional property that at time zero it is the zero section of $U$.

where in the last equation we used that $d_{(0,y)}\varphi$ maps the $i$-th canonical basis vector of $T_y\mathbb{R}^n \subset T_yU$ to $\vec{\boldsymbol{\alpha}}_i|_y$. Hence

$$\frac{d}{d\lambda}|_{\lambda=0}(\Phi \circ \mathbf{b}_\lambda) = \frac{d}{d\lambda}|_{\lambda=0}(\varphi \circ \mathbf{u}_\lambda) = \sum_i g_i \boldsymbol{\alpha}_i \tag{A.7}$$

lies in $\mathcal{B}|_V$, since $\boldsymbol{\alpha}_i \in \mathcal{B}|_V$ for all $i$. Now, any element of $\mathcal{S}$ is of the form $\frac{d}{d\lambda}|_{\lambda=0}(\Phi \circ \mathbf{b}_\lambda)$, for a 1-parameter family of global bisections $\mathbf{b}_\lambda$ s.t. $\mathbf{b}_0 = Id_M$. By Prop. A.5 we may assume that the $\mathbf{b}_\lambda$ have compact support. Equation (A.7) shows that there are finitely many open sets $V_i$ covering this support such that $\boldsymbol{\alpha}|_{V_i}$ lies in $\mathcal{B}|_{V_i}$. Therefore $\boldsymbol{\alpha} \in \mathcal{B}$, showing that $\mathcal{S} \subset \mathcal{B}$.

Having showed that $\mathcal{S} = \mathcal{B}$, we apply again Proposition 2.9 to conclude that Assumption 3.1 is satisfied. □

## A.4 The functoriality of differentiation

In this subsection we cast differentiation as a functor. We fix a manifold $M$ and consider two categories. They both differ slightly from those introduced in [31, §4].

The category $\mathsf{SingSub}_M$ has:

- objects:

$$\{(A, \mathcal{B}) \mid A \text{ a Lie algebroid over } M, \mathcal{B} \text{ a singular subalgebroid of } A\}$$

- arrows from $(A_1, \mathcal{B}_1)$ to $(A_2, \mathcal{B}_2)$:

$$\{\psi\colon A_1 \to A_2 \text{ a morphism of Lie algebroids covering } Id_M, \text{ such that } \psi(\mathcal{B}_1) \subset \mathcal{B}_2\}$$

The category $\mathsf{DiffeoGrd}_M$ has:

- objects:

$$\begin{aligned}\{\Phi\colon H \to G \text{ a morphism of diffeological groupoids covering } Id_M,\\ \text{where } H \text{ is diffeological groupoid which is holonomy-like,}\\ G \text{ is a source-connected Lie groupoid,}\\ \text{satisfying Assumption 3.1.}\}\end{aligned}$$

- arrows from $(\Phi_1\colon H_1 \to G_1)$ to $(\Phi_2\colon H_2 \to G_2)$:

$$\begin{aligned}\{(\Xi, F) \mid \Xi\colon H_1 \to H_2 \text{ a morphism of diffeological groupoids over } Id_M,\\ F\colon G_1 \to G_2 \text{ a morphism of Lie groupoids over } Id_M,\\ \text{s.t. the diagram below commutes}\}\end{aligned}$$

$$\begin{array}{ccc} H_1 & \xrightarrow{\Xi} & H_2 \\ {\scriptstyle \Phi_1}\downarrow & & \downarrow{\scriptstyle \Phi_2} \\ G_1 & \xrightarrow{F} & G_2 \end{array} \tag{A.8}$$

**Remark A.8.** One may view an object $(A, \mathcal{B})$ of the $\mathsf{SingSub}_M$ as the inclusion map $\mathcal{B} \hookrightarrow \Gamma_c(A)$. This point of view makes more manifest the parallel between the two above categories.

**Definition A.9.** We define the **differentiation functor**

$$\begin{aligned} D\colon \mathsf{DiffeoGrd}_M &\to \mathsf{SingSub}_M \\ (\Phi\colon H \to G) &\mapsto (A, \mathcal{B}) := (AG, \mathcal{S}) \\ (\Xi, F) &\mapsto F_*\colon AG_1 \to AG_2 \end{aligned}$$

where $\mathcal{S}$ was defined in §3.1, $AG$ denotes the Lie algebroid of $G$, and $(\Xi, F)$ is an arrow from $(\Phi_1\colon H_1 \to G_1)$ to $(\Phi_2\colon H_2 \to G_2)$.

The functor $D$ is well-defined on objects because of Thm. 3.2. We now check that $D$ is well-defined on arrows, i.e. that $F_*\mathcal{B}_1 \subset \mathcal{B}_2$. Every element of $\mathcal{B}_1$ is of the form $\frac{d}{d\lambda}|_{\lambda=0}(\Phi_1 \circ \mathbf{b}_\lambda)$ where $\{\mathbf{b}_\lambda\}_{\lambda\in I}$ is a 1-parameter family of global bisections for $H^{G_1}(\mathcal{B}_1)$ with $\mathbf{b}_0 = Id_M$. Using the commutativity of diagram (A.8) we have

$$F_*\left(\frac{d}{d\lambda}|_{\lambda=0}(\Phi_1 \circ \mathbf{b}_\lambda)\right) = \frac{d}{d\lambda}|_{\lambda=0}(F \circ \Phi_1 \circ \mathbf{b}_\lambda) = \frac{d}{d\lambda}|_{\lambda=0}(\Phi_2 \circ \Xi \circ \mathbf{b}_\lambda)$$

and the latter lies in $\mathcal{B}_2$, since $\Xi \circ \mathbf{b}_\lambda$ is a global bisection for $H^{G_2}(\mathcal{B}_2)$.

We now present some remarks and examples. At first sight the differentiation functor at the level of arrows looks strange: it maps $(\Xi, F)$ to $F_*$, hence $\Xi$ seems to "be lost". The following lemma shows that – under reasonable assumptions – this is not the case, since $\Xi$ is determined by $F$. Both assumptions appearing below are satisfied for integrals of singular subalgebroids, by Prop. 2.17.

**Lemma A.10.** Let $\Phi_i\colon H_i \to G_i$ be an object in $\mathsf{DiffeoGrd}_M$, for $i = 1, 2$. Assume that $H_1$ is source-submersive (Def. 2.15) and that $\Phi_2$ is coarsely injective (Def. 4.1 ). Let $(\Xi, F)$ be an arrow in $\mathsf{DiffeoGrd}_M$ from $\Phi_1$ to $\Phi_2$, as in diagram A.8. Then $\Xi$ is determined uniquely by $F$.

*Proof.* Assume that $(\Xi, F)$ and $(\tilde{\Xi}, F)$ are both morphisms in $\mathsf{DiffeoGrd}_M$ from $\Phi_1$ to $\Phi_2$. We have to show that $\Xi = \tilde{\Xi}$.

Take $\breve{H}_2$ to be a symmetric neighborhood of the identities in $H_2$ such that $\breve{H}_2 \cdot \breve{H}_2 \subset \overset{\circ}{H}_2$, where the latter is the neighborhood of the identities in $H_2$ satisfying the coarsely injective condition in Def. 4.1. Then $\overset{\circ}{H}_1 := \Xi^{-1}(\breve{H}_2) \cap \tilde{\Xi}^{-1}(\breve{H}_2)$ is a neighborhood of the identities in $H_1$.

Let $\mathbf{b}$ be a bisection in $\overset{\circ}{H}_1$. Then

$$(\Phi_2 \circ \Xi)\mathbf{b} = (F \circ \Phi_1)\mathbf{b} = (\Phi_2 \circ \tilde{\Xi})\mathbf{b}$$

due to the commutativity of diagram (A.8). This means that $\Xi\ \mathbf{b} * (\tilde{\Xi}\ \mathbf{b})^{-1}$ is a bisection in $\overset{\circ}{H}_2$ that maps under $\Phi_2$ to the identity bisection of $G_2$. The coarsely injective assumption on $\Phi_2$ implies that $\Xi\ \mathbf{b} = \tilde{\Xi}\ \mathbf{b}$. By the definition of source-submersive, every point of $H_1$ lies on a bisection, hence $\Xi$ and $\tilde{\Xi}$ agree on $\overset{\circ}{H}_1$. Finally, by the source-connectedness of $H_1$, they agree on the whole of $H_1$. □

In this remark we apply the functor $D$ to a class of examples.

**Remark A.11.** As we saw in Thm. A.7, given an object $(A, \mathcal{B})$ in $\mathsf{SingSub}_M$ such that $A$ is an integrable Lie algebroid, the holonomy groupoid $H^G(\mathcal{B})$ (for any choice of Lie groupoid $G$ of $A$) differentiates to $(A, \mathcal{B})$. We now extend this statement to arrows.

Consider objects $(A_1, \mathcal{B}_1)$ and $(A_2, \mathcal{B}_2)$ in $\mathsf{SingSub}_M$ together with an arrow between them, i.e. a Lie algebroid morphism $\psi\colon A_1 \to A_2$ such that $\psi_*(\mathcal{B}_1) \subset \mathcal{B}_2$. Make a choice of Lie groupoid morphism $F\colon G_1 \to G_2$ integrating $\psi$. (These data comprise exactly the arrows of the category $\mathsf{SingSub}_M^{Gpd}$ we defined in [31, §4]). We saw in [31, Thm. 4.6] that these data give rise canonically to an arrow $(\Xi, F)$ in $\mathsf{DiffeoGrd}_M$, where $\Xi\colon H^{G_1}(\mathcal{B}_1) \to H^{G_2}(\mathcal{B}_2)$. The arrow $(\Xi, F)$ differentiates to the original arrow $\psi\colon A_1 \to A_2$, by the definition of the functor $D$.

Loosely speaking, given a singular subalgebroid $\mathcal{B}$ of $AG$, the canonical morphism $\Phi\colon H^G(\mathcal{B}) \to G$ differentiates to the inclusion of $\mathcal{B}$ in $\Gamma_c(A)$. We make this precise as follows.

**Example A.12.** Let $\mathcal{B}$ be a singular subalgebroid of $A$, and $G$ a Lie groupoid integrating $A$. Consider the arrow $(\Phi, Id_G)$ in $\mathsf{DiffeoGrd}_M$ from $\Phi\colon H^G(\mathcal{B}) \to G$ to $Id_G$. It clearly differentiates to $Id_A\colon A \to A$, which at the level of sections restricts to the inclusion $\mathcal{B} \to \Gamma_c(A)$.

Notice that the arrow $Id_A$ in $\mathsf{SingSub}_M$ from $(A, \mathcal{B})$ to $(A, \Gamma_c(A))$, by the procedure recalled in Remark A.11, gives rise exactly to the arrow $(\Phi, Id_G)$ in $\mathsf{DiffeoGrd}_M$ considered in Ex. A.12. This was checked in [31, Ex. 4.12].

# B Rephrasing Assumption 3.1 as a liftability condition

In Def. 4.1 about integrals of singular subalgebroids, the differentiation condition a) contains in particular Assumption 3.1. In this appendix, for integrals, we give an equivalent characterization of this assumption in terms of the liftability of 1-parameter groups of global bisections of the Lie groupoid $G$, see Corollary B.2.

The following lemma addresses when bisections of $G$ can be "lifted" via $\Psi$.

**Lemma B.1.** Let $G$ be a Lie groupoid, $H$ be a diffeological groupoid, and $\Psi\colon H \to G$ a smooth morphism of diffeological groupoids covering $Id_M$, satisfying the coarsely injective condition c) in Def. 4.1.

i) There exists a neighborhood $\breve{H}$ of the identities in $H$ with the following property: for all bisections $\mathbf{c}$ of $G$, if there exists a bisection of $\breve{H}$ mapping to $\mathbf{c}$ under $\Psi$, then this bisection is unique.

ii) Let $\boldsymbol{\alpha} \in \Gamma_c(AG)$, consider the corresponding 1-parameter group $\mathbf{c}_\lambda := \exp_{1_M}(\lambda\overrightarrow{\boldsymbol{\alpha}})$ of global bisections of $G$. Assume that there is a 1-parameter family $\{\mathbf{b}_\lambda\}_{\lambda\in I}$ of global bisections of $H$ with $\mathbf{b}_0 = Id_M$ and mapping to $\{\mathbf{c}_\lambda\}_{\lambda\in I}$ under $\Psi$, where $I$ is an open interval containing 0. Then $\{\mathbf{b}_\lambda\}_{\lambda\in I}$ is a 1-parameter *group*, shrinking $I$ to a smaller open interval if necessary.

iii) For all $\boldsymbol{\alpha} \in \Gamma_c(AG)$, the following two conditions are equivalent:

   a) There is a 1-parameter *group* $\{\mathbf{b}_\lambda\}_{\lambda\in I'}$ of global bisections of $H$ such that $\frac{d}{d\lambda}|_{\lambda=0}(\Psi\circ\mathbf{b}_\lambda) = \boldsymbol{\alpha}$, where $I'$ is an open interval containing 0.

b) There is a 1-parameter family $\{\mathbf{b}_\lambda\}_{\lambda\in I}$ of global bisections of $H$ with $\mathbf{b}_0 = 1_M$ which maps to $\{\exp_{1_M}(\lambda\vec{\boldsymbol{\alpha}})\}_{\lambda\in I}$ under $\Psi$, where $I$ is an open interval containing 0.

*Proof.* i) Take $\breve{H}$ to be a symmetric neighborhood of the identity such that $\breve{H}\cdot\breve{H}\subset\overset{\circ}{H}$, where $\overset{\circ}{H}$ is as in Def. 4.1. Let $\mathbf{b}$ and $\mathbf{b}'$ be bisections of $\breve{H}$ with $\Psi(\mathbf{b}) = \mathbf{c} = \Psi(\mathbf{b}')$. Then $\mathbf{b} * (\mathbf{b}')^{-1}$ is contained in $\overset{\circ}{H}$, and maps to the identity bisection of $G$ under $\Psi$. Thus property c) in Def. 4.1 implies that $\mathbf{b} = \mathbf{b}'$.

ii) When $\lambda, \mu, \lambda+\mu \in I$ we compute

$$\Psi(\mathbf{b}_\lambda * \mathbf{b}_\mu) = \Psi(\mathbf{b}_\lambda) * \Psi(\mathbf{b}_\mu) = \mathbf{c}_\lambda * \mathbf{c}_\mu = \mathbf{c}_{\lambda+\mu} = \Psi(\mathbf{b}_{\lambda+\mu}).$$

Shrinking $I$ if necessary, we can assure that $\mathbf{b}_\lambda * \mathbf{b}_\mu$ and $\mathbf{b}_{\lambda+\mu}$ lie in $\breve{H}$. Item i) then implies that these two global bisections agree. Together with $\mathbf{b}_0 = 1_M$, this shows that $\{\mathbf{b}_\lambda\}_{\lambda\in I}$ is a 1-parameter group.

iii) Assume a). Since $\{\mathbf{b}_\lambda\}_{\lambda\in I}$ is a 1-parameter group and $\Psi$ is a groupoid morphism, $\{\Psi\circ\mathbf{b}_\lambda\}_{\lambda\in I}$ is also a 1-parameter group. Hence it agrees with $\{\exp_{1_M}(\lambda\vec{\boldsymbol{\alpha}})\}_{\lambda\in I}$, by the uniqueness of 1-parameter groups of bisections with initial velocity on a Lie groupoid.

Conversely, assume b). Using item ii), we see that a) follows, for some open subinterval $I' \subset I$. □

The following corollary[16] is an immediate consequence of Lemma B.1 iii).

**Corollary B.2.** In the Def. 4.1 of integral, the requirement that Assumption 3.1 is satisfied can be equivalently replaced by Assumption B.3 below.

**Assumption B.3.** For every $x \in M$ there is a neighborhood $V$ with this property:
For all $\boldsymbol{\alpha} \in \mathcal{B}$ with $Supp(\boldsymbol{\alpha}) \subset V$, there is a 1-parameter family $\{\mathbf{b}_\lambda\}_{\lambda\in I}$ of global bisections of $H$ with $\mathbf{b}_0 = Id_M$ which maps to $\{\exp_{1_M}(\lambda\vec{\boldsymbol{\alpha}})\}_{\lambda\in I}$ under $\Psi$.

[16] Notice that the openness condition b) in Def. 4.1 is not used to obtain this corollary.